\documentclass[a4paper,12pt,twoside]{article}
\usepackage{graphicx}\usepackage{a4wide}
\usepackage[pdftex,colorlinks,bookmarksnumbered,bookmarksopen,citecolor=red,urlcolor=red]{hyperref}
\usepackage{amssymb,amsmath,amsthm}
\usepackage{subfigure}%,color,colordvi,graphicx,xcolor}
\usepackage[only,llbracket,rrbracket,interleave]{stmaryrd}
\usepackage[sort,compress]{cite}
\usepackage{multirow}

\usepackage[makeroom]{cancel}

\newcommand{\bm}[1]{\text{\boldmath $#1$\unboldmath}}
\newcommand{\vect}[1]{\mathbf{#1}}
\newcommand{\mat}[1]{\mathbf{#1}}
\newcommand{\NN}{\mathbb{N}}

\newcommand{\sobo}[1][1]{\ensuremath{\mathcal{H}^{#1}}}
\newcommand{\eltwo}{\ensuremath{\mathcal{L}_2}}
\newcommand{\trial}{\ensuremath{\mathcal{V}}}
\newcommand{\test}{\sobo_{\Gamma^{\bmu}_D}}
\newcommand{\Dsp}[2][]{\bigl(#2\bigr)_{\! #1}}
\newcommand{\Bsp}[2][]{\bigl\langle#2\bigr\rangle_{#1}}
\newcommand{\nsd}    {\texttt{n}_{\texttt{sd}}}
\newcommand{\nnodes}{\texttt{n}_{\texttt{mn}}}
\newcommand{\npar}{\texttt{n}_{\texttt{pa}}}

\DeclareMathOperator{\adj}{adj}
\DeclareMathOperator{\tr}{tr}

\newcommand{\grad}{\bm{\nabla}}
\newcommand{\bn}{\bm{n}}
\newcommand{\bx}{\bm{x}}
\newcommand{\bX}{\bm{X}}

\newcommand{\bmu}{\bm{\mu}}
\newcommand{\bu}{\bm{u}}
\newcommand{\bv}{\bm{v}}

\newcommand{\bI}{\bm{\mathcal{I}}}
\newcommand{\I}{\mathcal{I}}
\newcommand{\dx}{\delta x}
\newcommand{\dy}{\delta y}
\newcommand{\dz}{\delta z}

\newcommand{\bd}{\bm{d}}

\newcommand{\A}{\mat{A}}
\newcommand{\0}{\mat{0}}
\newcommand{\stress}{\bm{\sigma}}
\newcommand{\deform}{\bm{\varepsilon}} 

\newcommand{\curve}[1][]{\bm{C}_{\!#1}} 
\newcommand{\ncp}    {\texttt{n}_{\texttt{cp}}}
\newcommand{\nk}      {\texttt{n}_{\texttt{k}}}
\newcommand{\surface}[1][]{\bm{S}_{\!#1}} 
\newcommand{\mcp}    {\texttt{m}_{\texttt{cp}}}

\newcommand{\nDeg}{\ensuremath{k}}

\newcommand{\Rin} {R_{\text{in}}}
\newcommand{\Rout}{R_{\text{out}}}
\newcommand{\Vin} {\Omega_{\text{in}}}
\newcommand{\Vout}{\Omega_{\text{out}}}

\usepackage[usenames]{color} %used for font color
\definecolor{clover}{RGB}{0,128,0}
\definecolor{purple}{RGB}{128,0,128}

\renewcommand{\emph}[1]{\textit{#1}}
\newenvironment{keywords}{\begin{quote}\emph{\textbf{Keywords:}}}{\end{quote}}
\newtheorem{remark}{Remark}

\begin{document}
%\begin{frontmatter}
%______________________________________________________________________
\title{Parametric solutions involving geometry integrated with computer-aided design}
%\author[ZCCE]{Ruben Sevilla}
%\ead{r.sevilla@swansea.ac.uk}
%\author[LaCaN]{Sergio Zlotnik}
%\ead{sergio.zlotnik@upc.edu}
%\author[LaCaN]{Antonio Huerta\corref{cor1}}
%\ead{antonio.huerta@upc.edu}
%\address[ZCCE]{Zienkiewicz Centre for Computational Engineering, College of Engineering,Swansea University, Bay Campus, SA1 8EN, Wales, United Kingdom.}
%\address[LaCaN]{Laboratori de Calcul Numeric (LaC\`aN), ETS de Ingenieros de Caminos, Canales y Puertos, Universitat Polit\`ecnica de Catalunya Barcelona, Spain.}
%\ead[url]{www.lacan.upc.edu}
\author{Ruben Sevilla\\[-1ex]
             \small Zienkiewicz Centre for Computational Engineering, \\[-1ex]
             \small College of Engineering, Swansea University, Bay Campus, SA1 8EN, Wales, UK \\[1em]
             Sergio Zlotnik, and Antonio Huerta\\[-1ex]
             \small Laboratori de C\`alcul Num\`eric (LaC\`aN), \\[-1ex]
             \small ETS de Ingenieros de Caminos, Canales y Puertos, \\[-1ex]
             \small Universitat Polit\`ecnica de Catalunya, Barcelona, Spain}

\date{\today}
%________________________________________________________________________
\maketitle

%______________________________________________________________________
\begin{abstract}

The main objective of this work is to describe a general and original approach for computing an off-line solution for a set of parameters describing the geometry of the domain. That is, a solution able to include information for different geometrical parameter values and also allowing to compute readily the sensitivities. Instead of problem dependent approaches, a general framework is presented for standard engineering environments where the geometry is defined by means of NURBS. The parameters controlling the geometry are now the control points characterising the NURBS curves or surfaces. The approach proposed here, valid for 2D and 3D scenarios, allows a seamless integration with CAD preprocessors. The proper generalised decomposition (PGD), which is applied here to compute explicit geometrically parametrised solutions, circumvents the curse of dimensionality. Moreover, optimal convergence rates are shown for PGD approximations of incompressible flows.
\end{abstract}
%________________________________________________________________________
%______________________________________________________________________
%\end{frontmatter}
%______________________________________________________________________
\section{Introduction}
	
The current role of computational simulations in modern engineering design is limited by the complexity of the simulations that are required, particularly during the final stages of a design. The main problem is motivated by the number of configurations that need to be tested (e.g. loads, boundary conditions, material parameters and geometric configurations). 
	
One alternative to decrease the computational complexity in this scenario is to introduce a reduced order model \cite{quarteroni2014reduced}. The main idea involves projecting the governing equations describing the full model onto a space with lower dimension that is described using a reduced order basis. Well known methods to produce reduced order basis are Krylov-based methods \cite{freund2003model}, the reduced basis method \cite{rozza2008reduced} and the proper orthogonal decomposition (POD) \cite{berkooz1993proper,lieu2006reduced}. More recently, the proper generalised decomposition (PGD) \cite{PGD-CCH:14,Chinesta-Keunings-Leygue} has gained popularity due to its ability to build reduced basis with no prior knowledge of the solution. The PGD starts by considering the solution not only as a function of the standard coordinates (i.e., space and time) but also of any parameter of interest (e.g.\ boundary conditions, external loads, material parameters). The problem involving a range of all the desired parameters can be solved at the cost of several problems of the same size as the original problem for a particular choice of the parameters. This expensive calculation, usually referred as the \textit{off-line stage}, is performed only once, usually making use of high performance computing resources, to build the reduced order basis that is written as an explicit function of the  coordinates (space and time) and the parameters (i.e.\ a computational vademecum, see \cite{FC-CLBACGAAH:13}). Then, the \textit{on-line stage} consists of just a particularisation of the solution by using, in the simplest case, an interpolation of the already computed results. 
	
The PGD has been successfully applied to numerous multi-dimensional problems involving boundary conditions, material parameters and external loads as extra coordinates, to name a few, see \cite{PGD-CCH:14,Chinesta-Keunings-Leygue} and references therein. The application to problems involving geometrically parametrised domains is generally more challenging and the existing work is usually limited to simple geometries where the extra coordinates are scalings: the length of an interval in one-dimensional problems \cite{Chinesta-CLBACGAAH:13}, the thickness of extruded geometries \cite{leygue2010first,bognet2012advanced} or problem specific parameters \cite{heuze2015parametric}. More recently, an approach based on an initial subdivision of the computational domain in macro-elements was proposed in \cite{AH-AHCCL:14}. This idea was extended to domains with interfaces in \cite{SZ-ZDMH:15} and has been also applied to an engineering design process in \cite{Neron-CNLB:16}. 
	
In this paper, a new approach to incorporate the geometric parameters as extra coordinates in a PGD framework is proposed. The objective is to produce a general methodology that enables to obtain the solution of a particular problem when the geometry of the domain is parametrised using the NURBS boundary representation of the domain, as usually done in a computer-aided design (CAD) environment. The control points of the NURBS entities defining the boundary are considered as extra coordinates and a mapping between a reference domain and the current configuration is proposed by using a solid mechanics analogy. An explicit and separated representation of the mapping is derived and the application of the proposed methodology to Stokes flow problems is presented using examples of increasing difficulty in two and three dimensions. Contrary to the approach in \cite{AH-AHCCL:14,SZ-ZDMH:15}, the methodology proposed in this paper considers geometric parameters that are independent on the spatial discretisation, i.e. the control points of the NURBS entities defining the boundary representation of the domain. In addition, the technique presented here enables the solution of the multi-dimensional problem to be computed using high-order finite elements whereas the technique introduced in \cite{AH-AHCCL:14,SZ-ZDMH:15} requires an affine mapping between a reference element and the macro-elements that are used to parametrise the geometry.
	
The structure of the remainder of the paper is as follows. Section \ref{sc:problemStatement} presents the problem statement using the Poisson equation on a geometrically parametrised domain and summarises the application of the PGD. Section \ref{sc:generalisedSolutions} describes in detail the proposed technique to build a generalised solution assuming that a mapping between a reference configuration and the current one can be written in separated form. In Section \ref{sc:cadIntegration}, a methodology to integrate this approach in a CAD environment is presented. Considering the geometric parameters as the control points of the NURBS entities describing the boundary of the domain, the methodology to build a mapping that can be explicitly written in separated form is detailed. Finally, Section \ref{sc:examples} presents a series of numerical examples of increasing difficulty involving the solution of Stokes flow problems in two and three dimensions.
	
%______________________________________________________________________
\section{Problem statement and geometrically parametrised solutions} \label{sc:problemStatement}
\subsection{The Poisson equation on a parametrised domain}
	
The methodology proposed here can be directly extended to second order linear problems such as the Stokes problem that is studied in the examples. However, in order to simplify the presentation, the heat problem (Poisson) is presented in detail for a parametrised domain $\Omega^{\bmu} \subset \mathbb{R}^{\nsd}$ (with $\nsd$ number of spatial dimensions), whose boundary $\partial\Omega^{\bmu}$ is characterised by a set of geometric parameters $\bmu\in\bI\subset\mathbb{R}^{\npar}$ (with $\npar$ number of parameters characterising the geometry) and is partitioned into Dirichlet, $\Gamma^{\bmu}_D$, and Neumann, $\Gamma^{\bmu}_N$, frontiers such that $\overline{\partial\Omega}^{\bmu}=\overline{\Gamma}^{\bmu}_D\cup\overline{\Gamma}^{\bmu}_N$ and $\Gamma^{\bmu}_D\cap\Gamma^{\bmu}_N=\emptyset$. Note that the set $\bI\subset\mathbb{R}^{\npar}$, which characterises the admissible range for parameters $\bmu$, can be defined as the Cartesian combination of the range for each parameter, namely,  $\bI := \I_1\times\I_2\times\dotsb\times\I_{\npar}$ with $\mu_i\in\I_i$ for $i=1,\dotsc , \npar$.
	
For each set of parameter $\bmu$, the objective is to determine the parametric solution $u^{\bmu}(\bx)$, with $\bx\in\Omega^{\bmu}$, of the boundary value problem
\begin{equation}\label{eq:heatStrong}
  \left\{
  \begin{aligned}
	- \grad\cdot(\mat{K}\grad u^{\bmu}) &= s       && \text{in $\Omega^{\bmu}$,}\\
	u^{\bmu} &= u_D   &&\text{on $\Gamma^{\bmu}_D$,}\\
	\bn\cdot(\mat{K}\grad u^{\bmu}) &= t         &&\text{on $\Gamma^{\bmu}_N$,}
  \end{aligned}
  \right.
\end{equation}
where $\mat{K}$ is the thermal conductivity (symmetric and positive definite) matrix, $s$ is a source term, $u_D$ is the imposed temperature,  $t$ is the imposed heat normal flux and $\bn$ is the outward unit normal vector.
The standard variational form of the previous problem reads: find $u^{\bmu}\in\trial^{\bmu}$ for all $v\in\test$ such that
\begin{subequations} \label{eq:nsd-D}
\begin{equation}\label{eq:heatVariational}
  a(\bmu ; u^{\bmu},v)=\ell(\bmu ; v) ,
\end{equation}
where the space of trial functions is $\trial^{\bmu}:= \{w\in \sobo(\Omega^{\bmu}) : w=u_D \text{ on }\Gamma^{\bmu}_D\}$ and its corresponding test functions space is $\test:= \{w\in\sobo(\Omega^{\bmu}) : w=0\text{ on }\Gamma^{\bmu}_D\}$.
The parametric bilinear and linear forms $a(\bmu ; \cdot,\cdot)$ and $\ell(\bmu ; \cdot)$ are defined by
\begin{equation} \label{eq:Forms}
  a(\bmu ; u,v):=\Dsp[\Omega^{\bmu}]{\grad v , \mat{K}\grad u}
  \text{ and }\;
  \ell(\bmu ; v):= \Dsp[\Omega^{\bmu}]{v , s} + \Bsp[\Gamma^{\bmu}_N]{v, t} ,
\end{equation}
\end{subequations}
where 
\begin{equation*}
  \Dsp[\Omega]{u , v} = \int_{\Omega} u \, v \, d\Omega \text{, }\;\quad 
  \Dsp[\Omega]{\bu , \bv} = \int_{\Omega} \bu\cdot\bv \, d\Omega  \text{, and }\;\quad
  \Bsp[\Gamma]{u , v}= \int_{\Gamma}{u \, v} \, d\Gamma,
\end{equation*}
denote, respectively, the $\eltwo$ product of scalar/vector functions in $\Omega$ and its traces over $\Gamma$.
	
%______________________________________________________________________
\subsection{The multi-dimensional parametric problem}
	
There are different alternatives to obtain, for any given set of parameters $\bmu$, an approximation to the solution $u^{\bmu}(\bx)$ of problem \eqref{eq:heatStrong}. The obvious option of solving a new problem for every instance of $\bmu$ is feasible but too costly. A standard strategy to reduce the cost is to pre-compute off-line some representative samples of the parametric family of solutions (e.g. snapshots for reduced basis methods, principal components for POD). Then, any other instance is computed on-line with a small computational overhead. Here, the PGD is preferred because the off-line phase provides an explicit description of the parametric solution, i.e.\ a computational vademecum see \cite{FC-CLBACGAAH:13}. Thus, in spite of an off-line phase more involved, the on-line phase is a simple functional evaluation with a negligible computational overhead. 
	
In practice, this can be interpreted as taking $\bmu$ as additional independent variables (or parametric coordinates) instead of problem parameters. Hence, the unknown temperature field is not interpreted any more as a parametric solution, denoted as $u^{\bmu}$,  but it is seen now as a function in a larger dimensional space and it is written as $u(\bx,\bmu)$ with  $(\bx , \bmu ) \in\Omega^{\bmu}\times\bI$. 
Consequently, formally $u$ lies in a tensor product space, namely, $u\in\trial^{\bmu}\otimes\eltwo(\I_1)\otimes\eltwo(\I_2)\otimes\dotsb\otimes\eltwo(\I_{\npar})$. A standard weighted residuals approach, with integrals in $\Omega^{\bmu}\times\bI$ and the usual integration by parts only in $\Omega^{\bmu}$ produces a weak form in this multi-dimensional setup. Namely, find $u\in\trial^{\bmu}\otimes\eltwo(\I_1)\otimes\eltwo(\I_2)\otimes\dotsb\otimes\eltwo(\I_{\npar})$ such that
\begin{subequations} \label{eq:Multi-D}
\begin{equation} \label{eq:weakform}
  A(u,v)=L(v),\;\forall v\in\test\otimes\eltwo(I_1)\otimes\eltwo(I_2)\otimes\dotsb\otimes\eltwo(I_{\npar}) ,
\end{equation}
with the following definitions of the bilinear and linear forms 
\begin{equation} \label{eq:Forms2}
  \begin{split}
  A(u,v)&:= \int_{\I_1}\!\int_{\I_2}\dotsi\int_{\I_{\npar}} a(\bmu ; u,v) \, d\mu_{\npar}\dotsm d\mu_2\, d\mu_1 \;\text{ and }\\
  L(v)&:=  \int_{\I_1}\!\int_{\I_2}\dotsi\int_{\I_{\npar}} \ell(\bmu ; v)\, d\mu_{\npar}\dotsm d\mu_2\, d\mu_1. 
  \end{split}
\end{equation}
\end{subequations}
	
Obviously, the number of dimensions of the solution domain increases with the number of parameters. To circumvent the \emph{curse of  dimensionality}, the PGD approach \cite{Ammar-AMCR:06,Chinesta-CLBACGAAH:13,PGD-CCH:14,Chinesta-Keunings-Leygue} is employed here. This approach assumes a separable structure in the function that approximates $u$. Note that the tensor product space $\trial^{\bmu}\otimes\eltwo(\I_1)\otimes\eltwo(\I_2)\otimes\dotsb\otimes\eltwo(\I_{\npar})$ inherits the multi-dimensional complexity of the problem and, in principle, does not assume separability of the functions. 
	
Moreover, the solution of \eqref{eq:weakform} requires an affine parameter dependence of the different forms. This is standard in reduced order methods and it is very well discussed in \cite{Patera-Rozza:07,Rozza:14}. More precisely, it is required that the different forms are expressed (or at least well approximated) by the sum of products of parameter-dependent functions and parameter-independent operators, for instance
\begin{equation*}
  a(\bmu ; u,v) = \sum_{q=1}^Q \biggl(\prod_{i=1}^{\npar} \Theta^q_i(\mu_i)\biggr) a^q(u,v) .
\end{equation*}
Note that the forms $a^q(u,v)$ do not depend on the parameters (in particular, they are integrated over domains parameter independent). In fact, finding the affine parameter dependence of \eqref{eq:Forms} is a major concern in subsequent sections and will enable to obtain a separated approximation of the solution, namely
\begin{equation*}
u \approx u_{\texttt{PGD}}^n = \sum_{m=1}^n \psi^u_m(\bmu) w_m(\bx).
\end{equation*}

%______________________________________________________________________
\section{Separated spatial mapping to determine generalised solutions}  \label{sc:generalisedSolutions}
	
If the affine parameter dependence must be enforced, it is necessary to integrate in space over domains not depending on the parameters. Note that forms $a(\bmu ; u,v)$ and $\ell(\bmu ; v)$, see \eqref{eq:Forms}, are integrated on spacial parametrised domains (domains depending on parameters $\bmu$). As suggested in \cite{AH-AHCCL:14,SZ-ZDMH:15}, a mapping $\bm{\mathcal{M}}_{\bmu}$ is necessary (not sufficient) in order to have an affine parameter dependence. This mapping transforms a reference domain $\Omega$ into the geometrically parametrised (``deformed'') domain $\Omega^{\bmu}$, namely
\begin{equation} \label{eq:mappingDomain}
  \begin{aligned}
  \bm{\mathcal{M}}_{\bmu}\, :
  \Omega\times\bI &\longrightarrow \Omega^{\bmu} \\
  (\bX,\bmu) &\longmapsto      \bm{\mathcal{M}}_{\bmu}(\bX,\bmu) = \bx = \bX + \bd(\bX,\bmu).
  \end{aligned}
\end{equation}
	
As classically in computational mechanics, the reference configuration $\Omega$ is associated to a reference coordinate system denoted by $\bX$, whereas the distorted domain $\Omega^{\bmu}$ will be associated to the spatial description $\bx$. Following this analogy, a \emph{displacement} field $\bd$ is used to relate both coordinate systems.
	
The mapping $\bm{\mathcal{M}}_{\bmu}$ can be defined in an \emph{ad hoc} manner for each problem \cite{AH-AHCCL:14,Neron-CNLB:16} or by a more general strategy \cite{AH-AHCCL:14,SZ-ZDMH:15}, but always has to induce an affine parameter dependence of forms \eqref{eq:Forms}. The introduction of this mapping allows to rewrite \eqref{eq:Forms} as integrals over the reference computational domain $\Omega$ and its corresponding boundary $\partial\Omega$, partitioned into the Dirichlet, $\Gamma_D$, and Neumann, $\Gamma_N$, boundaries, all independent of parameter $\bm{\mu}$. Namely,
\newcommand{\gradX}{\bm{\nabla}\!\!_{\bX}}
\newcommand{\Jaco}{\mat{J}_{\!\bmu}}
\begin{equation} \label{eq:FormsRef}
  \begin{split}
  a(\bmu ; u,v) &=\Dsp[\Omega]{\Jaco^{-1}\gradX v , \det(\Jaco)\mat{K}\Jaco^{-1}\gradX u }
	=\Dsp[\Omega]{\gradX v , \mat{H}_{\bmu}\,\gradX u }
  \text{ and }\\
  \ell(\bmu ; v) &= \Dsp[\Omega]{v , \det(\Jaco)\, s} + \Bsp[\Gamma_N]{v, \det(\Jaco)\, t} ,
  \end{split}
\end{equation}
where $[\Jaco]_{ij} = [\partial x_j/\partial X_i] $ is the Jacobian matrix of the mapping $\bm{\mathcal{M}}_{\bmu}$ and, more importantly, it is the only element in the above equations that depends upon parameters $\bm{\mu}$. Moreover, in order to compact the notation a new matrix $\mat{H}_{\bmu}$ is introduced,
\begin{equation}\label{eq:H}
  \mat{H}_{\bmu} := \frac{\adj(\Jaco^T)\,\mat{K}\,\adj(\Jaco)}{\det(\Jaco)},
\end{equation}
where the definition of the \emph{adjoint} operator, $\adj(\bm{A}) =\det(\bm{A})\,\bm{A}^{-1}$, has been used.
	
Note that even for \emph{ad hoc} \cite{AH-AHCCL:14} or a piecewise linear \cite{AH-AHCCL:14,SZ-ZDMH:15} mappings, the affine parameter dependence of the different forms in \eqref{eq:FormsRef} must be now determined. 

\begin{remark}[Mapping for the Stokes problem]
Despite being more cumbersome, the same rationale can be applied for the Stokes problem without any extra conceptual challenge. The bilinear viscosity form is reproduced here for illustration purposes
\begin{equation*}
  a(\bmu;\bm{u},\bm{v}) 
  = \Dsp[\Omega^{\bmu}]{ \grad  \bm{v}, \underline{\underline{\mat{C}}}        :\grad   \bm{u}}
  = \Dsp[\Omega]            { \gradX\bm{v}, \underline{\underline{\mat{\hat{C}}}}^{\bmu}:\gradX\bm{u}}, 
\end{equation*}
where $\text{C}_{ijkl}=\nu\,\delta_{ik}\,\delta_{jl}$, 
$\hat{\text{C}}_{ijkl}^{\bmu}= \nu\det(\mat{J_{\bmu}})\,\delta_{ik}\,[\mat{J_{\bmu}^{-T}}]_{js}\,[\mat{J_{\bmu}^{-1}}]_{sl}$ and $\nu>0$ is the kinematic viscosity.
Note that $[\mat{A}]_{ij}$ denotes the component $ij$ of a matrix $\mat{A}$.
\end{remark}
	
%______________________________________________________________________
\subsection{Separated displacements}

The methodology proposed here allows to obtain an affine parameter dependence (namely, a separable expression for $ \mat{H}_{\bmu}$) quasi-analytically if the mapping can be written (or at least, well approximated) with a separated representation, that is as a sum of separated terms, namely
\begin{equation}\label{eq:MapSep}
  \bm{\mathcal{M}}_{\bmu}(\bX,\bmu) = \bX + \sum_{m=1}^N \psi_m(\bmu)\; \bd^m(\bX) ,
\end{equation}
which induces a separated Jacobian matrix 
\begin{equation}\label{eq:JacoSep}
%\begin{split}
  \Jaco(\bX,\bmu) = \frac{\partial\bx}{\partial\bX}(\bX,\bmu) 
%                  = \mat{I}_{\nsd} + \sum_{m=1}^N \psi_m(\bmu)\; \frac{\partial\bd^m}{\partial\bX}(\bX)
                  = \mat{I}_{\nsd} + \sum_{m=1}^N \psi_m\,\mat{A}_m ,
%\end{split}
\end{equation}
where $\mat{A}_m := [\partial\bd^m/\partial\bX]$ for $m=1,\dotsc , N$.
	
%Thus, if the mapping $\bm{\mathcal{M}}_{\bmu}$ is separated as shown in \eqref{eq:MapSep}, then the new matrix $\mat{H}_{\bmu}$, see \eqref{eq:H}, can be written in separated form and, consequently, the forms in \eqref{eq:FormsRef} will have an affine parameter dependence.

%______________________________________________________________________
\subsection{Affine parameter dependence}

The affine parameter dependence of $a(\bmu ; u,v)$ as defined in  \eqref{eq:FormsRef}, is in practice determined by obtaining a separated expression of $\mat{H}_{\bmu}$. 
In order to obtain a separated expression for matrix $\mat{H}_{\bmu}$, both $\adj(\Jaco)$ and $\det(\Jaco)$ are analytically separated. Then, the Higher-Order PGD-Projection proposed in \cite{DM-MZH:15} is used to obtain a compact separation for $\mat{H}_{\bmu}$.

The separated representation for the determinant, $\det(\Jaco)$, can be obtained using \emph{Leibniz formula} from equation \eqref{eq:JacoSep}, whereas for $\adj(\Jaco)$, the Leverrier's algorithm \cite{FaddeevBook} is employed. This method is a consequence of the Cayley-Hamilton theorem \cite{cayley1858memoir,Householder} and the Newton's identities \cite{mead1992newton}. It enables to express the adjoint of a matrix $\mat{A} \in \mathbb{R}^{n \times n}$ in terms of its trace and its powers, namely
\begin{equation}\label{eq:InvByC-H}
  \adj(\mat{A}) = \sum_{s=0}^{n-1} \mat{A}^s 
	\sum_{k_1,k_2,\dotsc ,k_{n-1}}\prod_{l=1}^{n-1}
	\frac{(-1)^{k_l+1}}{l^{k_l}k_l!}\tr(\mat{A}^l)^{k_l}
\end{equation}
where $k_l\in\NN_0$ and $s+\sum_{l=1}^{n-1}l k_l =n-1$. In practice, the two cases of interest are 2D ($n=2$) and 3D ($n=3$) problems, which are detailed next.
	
%______________________________________________________________________
\subsection{Two-dimensional approach}
	
In 2D the Jacobian $\Jaco$ is a $2\times 2$ matrix. From \eqref{eq:JacoSep}, the separated expression for its determinant is
\begin{multline}\label{eq:detJaco2D}
\det(\Jaco) = \Bigl(1+\sum_{m=1}^N \psi_m(\bmu) \bigl[\mat{A}_m\bigr]_{11} \Bigr)
                     \Bigl(1+\sum_{m=1}^N \psi_m(\bmu) \bigl[\mat{A}_m\bigr]_{22} \Bigr) \\- 
                     \Bigl(    \sum_{m=1}^N \psi_m(\bmu) \bigl[\mat{A}_m\bigr]_{21} \Bigr)
                     \Bigl(    \sum_{m=1}^N \psi_m(\bmu) \bigl[\mat{A}_m\bigr]_{12} \Bigr) ,
\end{multline}
again $[\mat{A}]_{ij}$ denotes the component $ij$ of the matrix $\mat{A}$.

To obtain the expression for the adjoint requires to particularise \eqref{eq:InvByC-H}, namely
\begin{equation*}
  \adj(\mat{A}) =  \tr(\mat{A})\mat{I}_2 - \mat{A} , 
\end{equation*}
which is a linear mapping because the trace is also linear. 
%
%\begin{equation*}
% \adj\Bigl(\sum_{i=1}^n\alpha_i\mat{A}_i\Bigr) = \sum_{i=1}^n\alpha_i\adj(\mat{A}_i) ,
%\end{equation*}
%
%because the trace of any square matrix is a linear mapping. 
Thus, the Jacobian in separated form as presented in \eqref{eq:JacoSep} can be written as 
\begin{equation}\label{eq:AdJaco2D}
  \adj(\Jaco) = \mat{I}_2 + \sum_{m=1}^N \psi_m\,\adj\bigl([\partial\bd^m/\partial\bX]\bigr) 
              = \mat{I}_2 + \sum_{m=1}^N \psi_m\,\adj(\mat{A}_m) .
\end{equation}
	
Consequently, recalling that $\adj(\mat{A}^T)=\adj(\mat{A})^T$, the matrix $\mat{H}_{\bmu}$ in \eqref{eq:H} can be rewritten as
\begin{multline*}%\label{eq:H}
  \det(\Jaco)\; \mat{H}_{\bmu} = \mat{K}  
%                                                            + \sum_{m=1}^N  \psi_m\,\adj\bigl([\partial\bd^m/\partial\bX]\bigr)^T \,\mat{K}
  + \sum_{m=1}^N \psi_m\,\mat{K}\,\adj(\mat{A}_m) 
  + \sum_{m=1}^N \psi_m\bigl[\mat{K}\,\adj(\mat{A}_m)\bigr]^T \\
  + \sum_{m=1}^N\sum_{l=1}^N\psi_m\psi_l\,\adj(\mat{A}_m)^T\,\mat{K}\,\adj(\mat{A}_l) .
\end{multline*}
The separated expression of $\mat{H}_{\bmu}$ is efficiently obtained with a numerical Higher-Order PGD-Projection \cite{DM-MZH:15}.

%______________________________________________________________________
\subsection{Three-dimensional approach}

Following the previous rationale, in 3D the Jacobian $\Jaco$ is a $3\times 3$ matrix. The separated expression for the determinant, $\det(\Jaco)$, is obtained using \emph{Leibniz formula} as 
\begin{equation}\label{eq:detJaco3D}
 \det(\Jaco) = \sum_{\sigma\in\mathcal{S}_3}\operatorname{sgn}(\sigma)
      \prod_{i=1}^3 \Bigl(1{+}\sum_{m=1}^N \psi_m(\bmu) \bigl[\mat{A}_m\bigr]_{i\sigma(i)} \Bigr),
\end{equation}
where $\mathcal{S}_3$ is the set of the six permutations of the integers $\{1, 2, 3\}$, where the element in position $i$ after the reordering $\sigma$ is denoted $\sigma(i)$, and $\operatorname{sgn}(\sigma)$ denotes the signature of $\sigma$ (i.e.\ $+1$ for even $\sigma$ and $-1$ for odd $\sigma$). Note that such a separation will induce sums of order $N^3$.

The adjoint of the Jacobian is also separated by means of particularising \eqref{eq:InvByC-H} to 3D,  
\begin{equation*}
\adj(\mat{A}) = \tfrac{1}{2} \bigl[\bigl(\tr(\mat{A})\bigr)^2-\tr(\mat{A}^2)\bigr]\mat{I}_3 - \tr(\mat{A})\mat{A} + \mat{A}^2 , 
\end{equation*}
and given the separation of the Jacobian in \eqref{eq:JacoSep} implies 
\begin{multline}\label{eq:AdJaco3D}
\adj(\Jaco) = \mat{I}_3 + \sum_{m=1}^{N}\psi_m \bigl[ \tr(\mat{A}_m)\mat{I}_3 - \mat{A}_m \bigr] \\
+ \sum_{m=1}^{N} \sum_{l=1}^{N} \psi_m \psi_l \bigl[ \frac{1}{2} \bigl( \tr(\mat{A}_m)\tr(\mat{A}_l) - \tr(\mat{A}_m \mat{A}_l) \bigr) \mat{I}_3 - \tr(\mat{A}_m) \mat{A}_l + \mat{A}_m \mat{A}_l  \bigr].
\end{multline}
As done in 2D, once the adjoint is separated, from \eqref{eq:H}, a separated expression for 
$\det(\Jaco)\mat{H}_{\bmu} = \adj(\Jaco^T)\,\mat{K}\,\adj(\Jaco)$ can be computed. In this case, it consists of $N^4$ terms. As previously noted for 2D, the separated expression of $\mat{H}_{\bmu}$ is obtained with the Higher-Order PGD-Projection described next.

%______________________________________________________________________
\subsection{The Higher-Order PGD-Projection} \label{sc:HOPGD}

Finally, to obtain a separated approximation of $\mat{H}_{\bmu}$, namely
\begin{equation*}
 \mat{H}_{\bmu}\approx\mat{H}_{\bmu}^{\texttt{sep}}(\bx,\bmu) 
                          = \sum_{m=1}^{N_H}  \psi^{\mat{H}}_{m}(\bmu)\, \mat{A}^{\mat{H}}_m(\bX),
\end{equation*}
equation \eqref{eq:H} has to be solved. This is done by means of the Higher-Order PGD-Projection \cite{DM-MZH:15}, which in practice obtains the separated approximation with an $\eltwo$ projection, say
\begin{equation*}
    \bigl(\bm{v},\mat{H}^{\texttt{sep}}_{\bmu}\bigr)_{\Omega\times\bI} 
 = \bigl(\bm{v},\adj(\Jaco^T)\,\mat{K}\,\adj(\Jaco)/\det(\Jaco)\bigr)_{\Omega\times\bI}  
\end{equation*}
for all $\bm{v}$ in a suitable space. In practice, the PGD strategy is implemented. Thus, a greedy approach is used. For computational efficiency, it is critical to use the exact separated expressions for the determinant, see \eqref{eq:detJaco2D} and \eqref{eq:detJaco3D}, and the adjoint of the Jacobian, see \eqref{eq:AdJaco2D} and \eqref{eq:AdJaco3D}. In practice, each mode of $\mat{H}_{\bmu}^{\texttt{sep}}$ is obtained by solving
\begin{multline} \label{eq:hprojection2}
    \bigl(\bm{v},\det(\Jaco) \, \psi^{\mat{H}}_{n}\, \mat{A}^{\mat{H}}_n\bigr)_{\Omega\times\bI} 
 = \bigl(\bm{v}, \adj(\Jaco^T)\,\mat{K}\,\adj(\Jaco)\bigr)_{\Omega\times\bI} \\
 - \bigl(\bm{v},\det(\Jaco) \, \sum_{m=1}^{n-1}  \psi^{\mat{H}}_{m}\, \mat{A}^{\mat{H}}_m\bigr)_{\Omega\times\bI},
\end{multline}
where now, as usual in PGD, $\bm{v}$ lives in the tangent space of the mode. Note that the Higher-Order PGD-Projection allows many parameters without any prior knowledge, is computationally efficient, and the precision of this approximation (i.e.\ the total number of terms $N_H$) can be controlled by the user.

%______________________________________________________________________
\section{Integration within a CAD environment}  \label{sc:cadIntegration}
	
In this section a procedure to integrate the methodology described in the previous section within a CAD environment is proposed. To simplify the presentation, the two dimensional case is presented here and the details for the extension to three dimensional domains are given in~\ref{app:cadIntegration3D}.
	
The boundary of the parametrised domain, $\partial \Omega^{\bmu}$, is assumed to be described by a set of NURBS curves $\{\curve[j]^\bmu\}_{j=1,\dotsc, M}$, being $M$ the total number of curves, namely
\begin{equation*}
  \partial \Omega^{\bmu} =  \bigcup_{j=1}^{M} \curve[j]^\bmu([0,1]).
\end{equation*}
	
Next, the necessary concepts about NURBS curves are briefly recalled and the proposed strategy to build a geometric mapping $\bm{\mathcal{M}}_{\bmu}$ that can be written in the separated form~\eqref{eq:MapSep} is presented in detail.
	
%______________________________________________________________________
\subsection{NURBS curves}  \label{sbc:nurbsCurves} 
	
A $q$th-degree non-uniform rational B-spline (NURBS) curve is a piecewise rational function defined in parametric form as
\begin{equation} \label{eq:nurbsDefinition}
  \curve(\lambda) = \sum_{i=0}^{\ncp} \bm{B}_i\,R_i(\lambda)
  \qquad \lambda \in [0,1]
  \end{equation}
where $\{\bm{B}_i\}$ are the coordinates of the $\ncp+1$ \emph{control points} (forming the \emph{control polygon}) and $\{R_i(\lambda)\}$ are rational basis functions defined as
\begin{equation*} \label{eq:rationalBasis}
  R_i(\lambda) =  \nu_i \, C_i^q(\lambda)
  \biggm/
  \left(\sum_{i=0}^{\ncp} \nu_i\,          C_i^q(\lambda)\right) .
\end{equation*}
	
In the above expression $\{\nu_i\}$  are the control weights associated to the control points and $\{C_i^q(\lambda)\}$ are the normalized B-spline basis functions of degree $q$, which are defined recursively by
\begin{align*}
  C_i^0(\lambda) &= \begin{cases}
	1 & \text{if $\lambda\in    [\lambda_i,\lambda_{i+1})$}  \\
	0 & \text{elsewhere}  \\
	\end{cases}\\
	C_i^k(\lambda) &= \frac{\lambda-\lambda_i}{\lambda_{i+k}-\lambda_i}
	C_i^{k-1}(\lambda)+
	\frac{\lambda_{i+k+1}-\lambda}{\lambda_{i+k+1}-\lambda_{i+1}}
	C_{i+1}^{k-1}(\lambda)
\end{align*}
for $k=1,\dotsc , q$ and where $\lambda_i$ (for $i=0,\dotsc, \nk$) are the \emph{knots} or \emph{breakpoints}, which are assumed ordered $0\le\lambda_i\le\lambda_{i+1}\le 1$. They form the so-called \emph{knot vector},
\begin{equation*}
  \Lambda = \{ \underbrace{0,\ldots,0}_{q+1},
  \lambda_{q+1},\dotsc,\lambda_{\nk-q-1},
  \underbrace{1,\dotsc ,1}_{q+1} \} ,
\end{equation*}
which uniquely describes the B-spline basis functions.  
%The multiplicity of a knot, when it is larger than one, determines the decrease in the number of continuous derivatives. 
The number of control points, $\ncp+1$, and knots, $\nk+1$, are related to the degree of the parametrisation, $q$, by the relation $\nk = \ncp + q + 1$, see \cite{Piegl-Tiller:95} for more details.
	
\begin{figure}[!tb]
	\centering
	\subfigure[]{\includegraphics[width=0.49\textwidth]{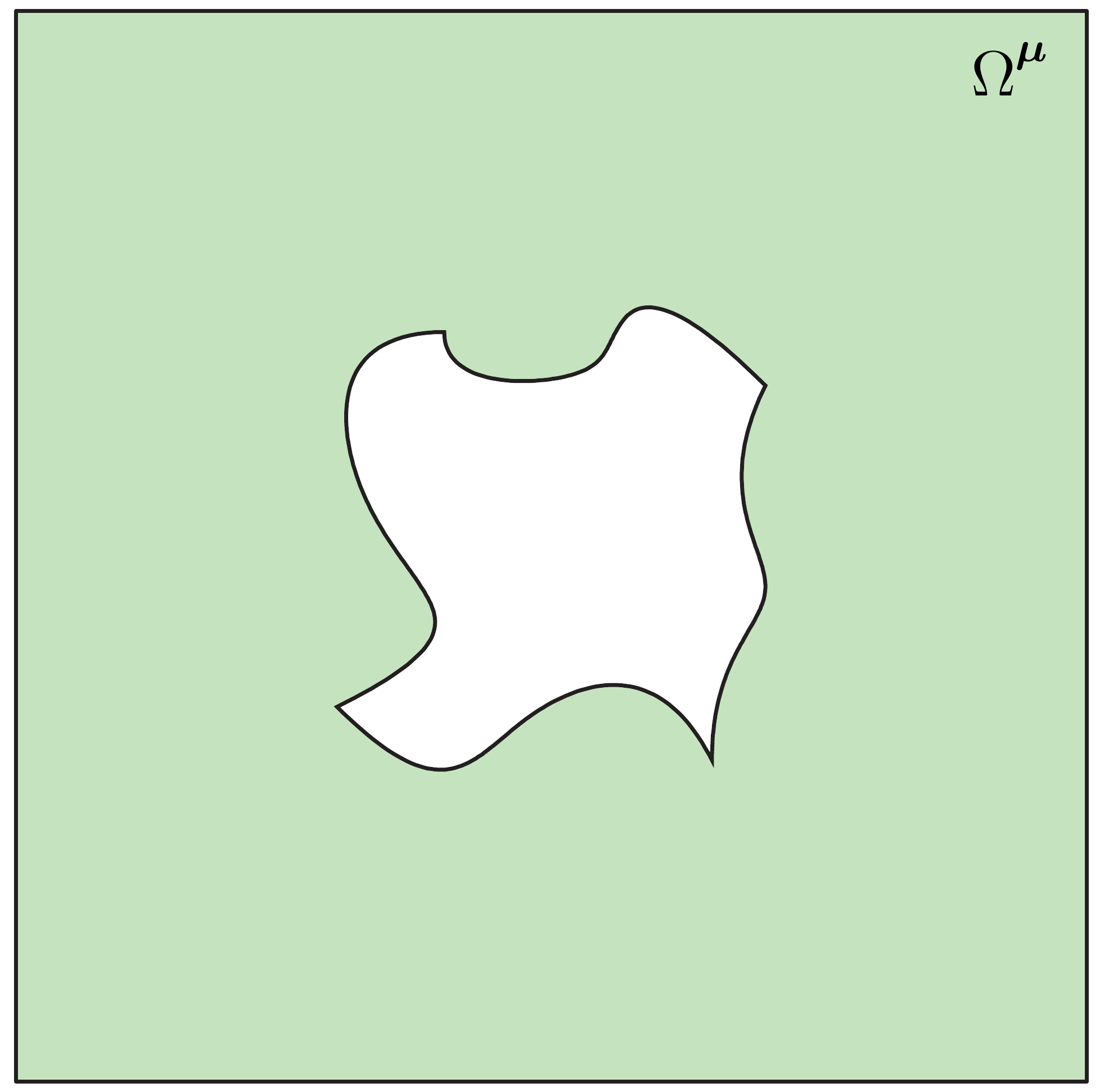}}
	\subfigure[]{\includegraphics[width=0.49\textwidth]{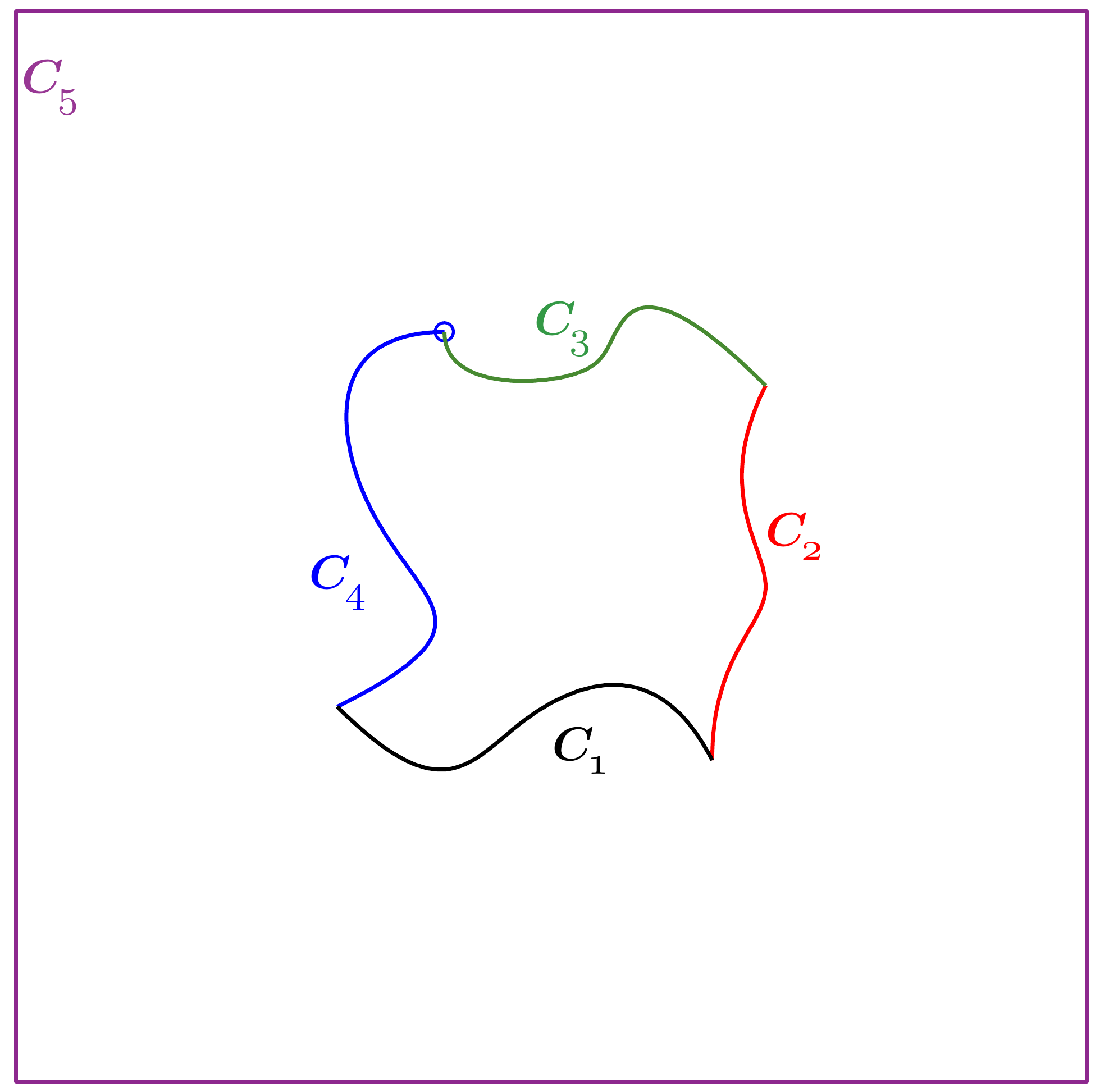}}
	\caption{(a) Domain $\Omega^{\bmu}$ and (b) NURBS curves describing the boundary of $\Omega^{\bmu}$, where each colour represents a different curve.}
	\label{fig:domainAndBoundary}
\end{figure}

\begin{figure}[!tb]
  \centering
  \includegraphics[width=0.75\textwidth]{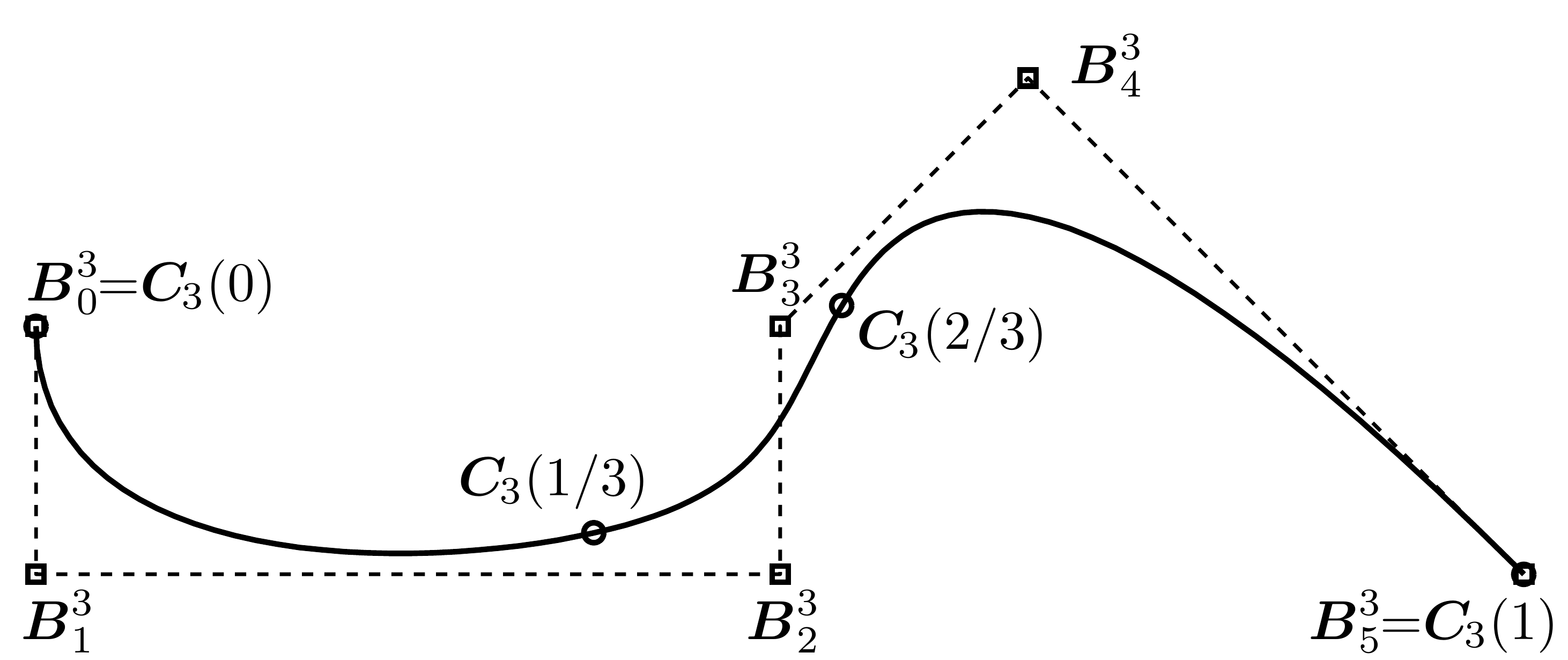}
  \caption{NURBS curve (solid line), control points (denoted by $\Box$), control polygon (dashed line) and breakpoints (denoted by $\circ$).}
  \label{fig:nurbsCurve}
\end{figure}

Figure \ref{fig:domainAndBoundary} shows an example of a two dimensional domain $\Omega^{\bmu}$ where the boundary is described by five NURBS curves.
The curve $\curve[3]$ in Figure \ref{fig:domainAndBoundary} (b) is represented in Figure \ref{fig:nurbsCurve} with the corresponding control polygon, formed by six control points and the breakpoints. The knot vector for this curve is given by 
\begin{equation*}
  \Lambda=\{0,0,0,1/3,2/3,1,1,1\}.
\end{equation*}
	
%______________________________________________________________________
\subsection{Geometric parameters}  \label{sbc:geometricParameters}
\newcommand{\dc}{\delta\!\bm{B}}%{\bmu}
	
The geometric parameters $\bmu\in\bI\subset\mathbb{R}^{\npar}$ are defined as the variations of the original coordinates of the control points of the NURBS curves describing the boundary. More precisely, for each NURBS curve $\curve[j]$, with $j=1,\dotsc, M$, having $\ncp^j + 1$  control points, the undisturbed boundary is characterised by the coordinates of the control points:  $\bm{B}_i^j$, for  $i=0,\dotsc , \ncp^j$. This configuration will be used as the reference one in $\Omega$ and will be associated to a reference coordinate system denoted by $\bX$. The distorted domain, $\Omega^\bmu$, will be associated to the spatial description $\bx$. The boundary in the spatial domain, $\partial\Omega^\bmu$, is defined by the position of the displaced control points, namely $\bm{B}_i^j+\dc_i^j$. The displacement range for each control point is characterised by 
\begin{equation*}
  \dc_i^j\in\bI^j_i = [\underline{\dx}_i^j, \overline{\dx}_i^j] \times [\underline{\dy}_i^j, \overline{\dy}_i^j].
\end{equation*}
In fact, each displacement of a control point $i$ on the $j$-th NURBS, $\dc_i^j$, might depend upon the parameters and can be written as 
\begin{equation}\label{eq:DeltaControl}
  \dc_i^j=\mu_1^{i,j}\bm{e}_1+\mu_2^{i,j}\bm{e}_2 ,
\end{equation}
where $\bm{e}_i$, for $i=1,2$, are the unit coordinate vectors. Then $\bmu\in\bI := \bI^1\times\bI^2 \times\dotsb\times\bI^M$, where $\bI^j = \bI^j_1 \times \bI^j_2 \times \cdots \bI^j_{\ncp^j + 1}$ is the range of variation of the coordinates of the control points of the curve $\curve[j]$. 
	
Consequently, the maximum number of geometric parameters is
\begin{equation} \label{eq:maxParams}
  \npar = \sum_{j=1}^M (\ncp^j + 1) \nsd ,
\end{equation}
but in practical problems, the number of geometric parameters $\npar$ is drastically lower than $m$ because not all the control points of all the curves are to be modified during the design stage.

\begin{remark} \label{rk:reparametrisation}
In a practical setting it is common to introduce some restrictions on the motion of the control points (viz.\ pure translations, rotations, expansions...), meaning that the motion of a set of control points can be expressed with a significantly low number of parameters.
\end{remark}
	
%______________________________________________________________________
\subsection{Separated representation of the boundary displacement}  \label{sbc:boundaryDisplacement}
	
The variation of a control point $\bm{B}_i^j$ of a NURBS curve $\curve[j]$, namely $\dc_i^j$, changes the definition of the original curve only in the support of the basis function $R_i^j$, given by the subspace of the parametric space $[\lambda_i, \lambda_{i+q^j+1}]$. The modified NURBS curve is parametrised by
\begin{equation} \label{eq:nurbsDefinitionChangeCP}
  \curve[j]^{\bmu}(\lambda) = \sum_{i=0}^{\ncp^j}  ( \bm{B}_i^j + \dc_i^j ) \,R_i^j(\lambda)
  \qquad \lambda \in [0,1] .
\end{equation}
	
\begin{figure}[!tb]
  \centering
  \includegraphics[width=0.75\textwidth]{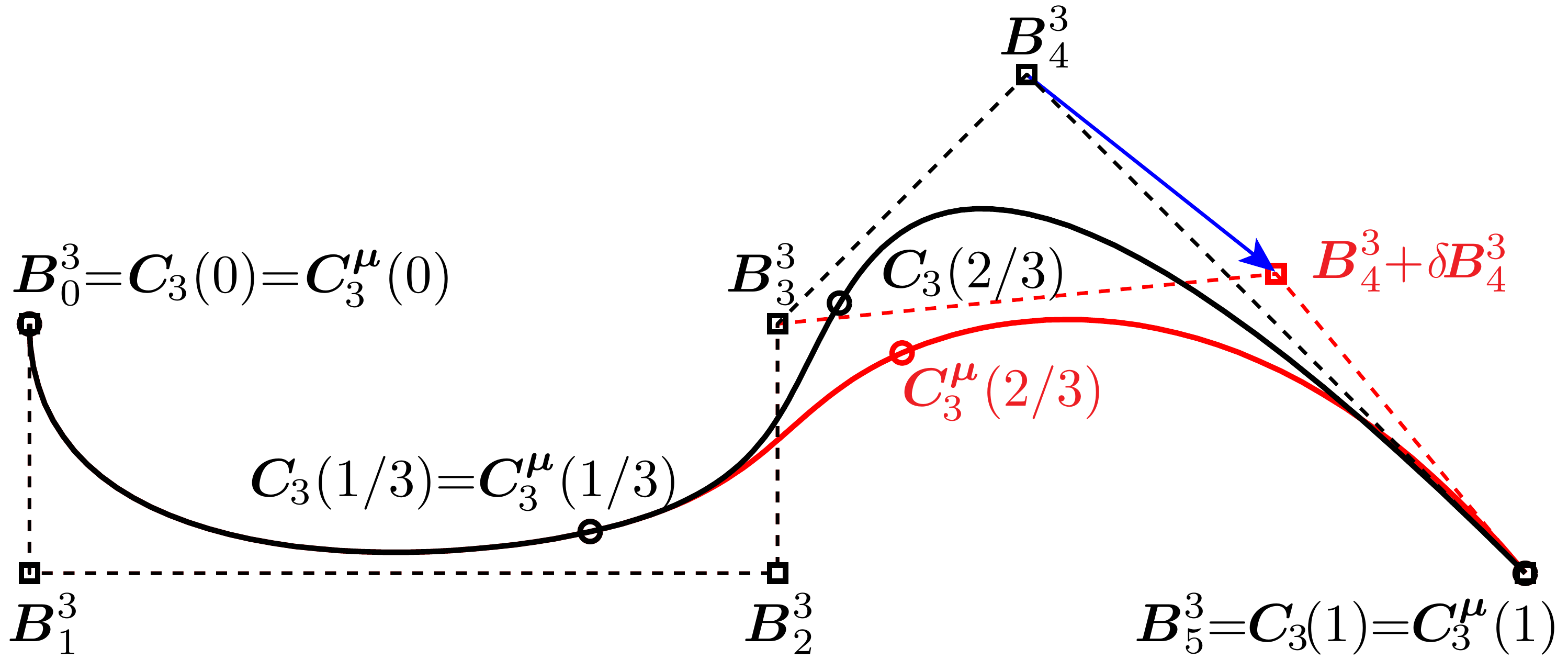}
  \caption{Change in the shape of a NURBS curve induced by the variation of the coordinates of the control point $\bm{B}_4^3$ given by $\dc_4^3$.}
  \label{fig:changeControlPoint}
\end{figure}

Figure \ref{fig:changeControlPoint} illustrates the effect of modifying the coordinates of one control point of a NURBS. The curve in red is the result of modifying the coordinates of the control point $\bm{B}_4^3$ of the original curve in black, also depicted in Figure \ref{fig:nurbsCurve}. It can be observed that both curves are identical in the interval [0,1/3] of the parametric space whereas they differ in the interval [1/3,1], which is the support of the basis function $R_4^3$ associated to the control point $\bm{B}_4^3$.
	
Given a computational mesh for the reference configuration, $\Omega$, the boundary mesh nodes affected by the motion of a control point can be easily identified. The procedure starts by finding the NURBS curve to which each boundary mesh node belongs and its associated parametric coordinate by using a standard NURBS point projection algorithm \cite{Piegl-Tiller:95}. For each boundary mesh node $\bX_k$, the index $j\in\{1, \dotsc ,M\}$ and parametric coordinate $\lambda_k$ such that $\curve[j](\lambda_k) = \bX_k$ are computed. Then, for a deformed configuration, induced by the variation of a control point $\bm{B}_i^j$ of a NURBS curve $\curve[j]$, namely $\dc_i^j$, the new position of each boundary node is computed as $\bx_k^\bmu = \curve[j]^\bmu(\lambda_k)$. 
	
Therefore, the variation of the control points of a boundary curve $\curve[j]$ induces a displacement of the boundary mesh nodes, namely 
\begin{equation*}
  \delta\bd^j(\bX_k ,\bmu) = \bx_k^\bmu-\bX_k \quad\text{, for all $k \in \mathcal{S}$,}
\end{equation*}
where $\mathcal{S} = \left\{ l \in \left\{1,2,\ldots, \nnodes \right\} : \bX_l \in \partial \Omega \right\}$ is the index set of mesh nodes on the boundary of the computational domain and $\nnodes$ is the total number of mesh nodes. Using the expression of the original and modified NURBS boundary curves, Equations \eqref{eq:nurbsDefinition} and \eqref{eq:nurbsDefinitionChangeCP} respectively, the displacement of the boundary mesh node $\bX_k=\curve[j](\lambda_k)$ that belongs to the NURBS curve $\curve[j]$ can be written in separated form as
\begin{equation*}\label{eq:boundayDisplacement}
  \delta\bd^j(\bX_k ,\bmu) = \sum_{i=0}^{\ncp^j} \dc_i^j \,R_i^j(\lambda_k) 
                                        = \sum_{i=0}^{\ncp^j} \sum_{s=1}^{\nsd} \mu_s^{i,j}\bm{e}_s\,R_i^j(\lambda_k) ,
\end{equation*}
where the dependence of the displacements of the control points in terms of the parameters described in \eqref{eq:DeltaControl} has been used. Moreover, since the NURBS parameter $\lambda_k$ is only dependent on the spatial coordinates $\bX_k$, and not on the geometric parameters $\bmu$, the previous equation can be written as,
\begin{equation*} %\label{eq:boundayDisplacementSeparated}
  \delta\bd(\bX_k ,\bmu) = \sum_{j=1}^{M} \delta\bd^j(\bX_k ,\bmu) 
       = \sum_{j=1}^{M} \sum_{i=0}^{\ncp^j}\sum_{s=1}^{\nsd} \mu_s^{i,j}\bm{e}_s\,R_i^j\bigl(\curve[j]^{-1}(\bX_k)\bigr) ,
\end{equation*}
which characterises the displacement of the boundary nodes and has a separated expression of the form
\begin{equation} \label{eq:boundayDisplacementSeparated}
  \delta\bd(\bX_k ,\bmu) = \sum_{i=1}^{\npar} \phi_i(\mu_i) \bd^{\texttt{b}}_i(\bX_k) .
\end{equation}

Note that this expression is compatible with the desired structure of the displacements described in \eqref{eq:MapSep}. In fact, it can be further compacted by means of the Higher-Order PGD-Projection \cite{DM-MZH:15} to obtain a more compact separation in the form of 
\begin{equation*}
  \delta\bd(\bX_k ,\bmu) = \sum_{i=1}^{\npar} \phi_i(\mu_i) \bd^{\texttt{b}}_i(\bX_k)
                                      = \sum_{m=1}^N\prod_{i=1}^{\npar} \widehat{\psi}_i^m(\mu_i) \bm{\hat{d}}^{\texttt{b}}_i(\bX_k)
                                      = \sum_{m=1}^N \psi_m^{\delta}(\bmu) \bm{\hat{d}}^{\texttt{b}}_i(\bX_k) ,
\end{equation*}
which now coincides with \eqref{eq:MapSep} but only for the boundary nodes. Precisely, the next section describes the extension to any point in the domain.

It is also important to note that the total number of terms in the separated representation \eqref{eq:boundayDisplacementSeparated} will be, in the majority of cases, much lower than $m$, defined in \eqref{eq:maxParams}. In practical applications, a large percentage of mesh boundary nodes will not be affected by a variation of the control points of a particular NURBS curve describing the boundary, so the displacement will be zero for a large number of boundary nodes. More precisely, the set of mesh nodes affected by the variation of a control point $\bm{B}_i^j$ of a NURBS curve $\curve[j]$ can be defined as
\begin{equation} \label{eq:boundaryNodesAffected}
  \mathcal{S}_i^j = \left\{ l \in \mathcal{S} : \bX_l = \curve[j](\lambda_l) \text{ for } \lambda_l \in [\lambda_i, \lambda_{i+q+1}] \right\}.
\end{equation}

%______________________________________________________________________
\subsection{Separated representation of the geometric mapping}   \label{sbc:domainDisplacement}
	
The proposed strategy to build a mapping $\bm{\mathcal{M}}_{\bmu}$ between the reference configuration, $\Omega$, and the current configuration, $\Omega^\bmu$, consists on solving a solid mechanics problem. The reference configuration is assumed to be a linear elastic medium and the displacement of mesh boundary nodes, induced by the variation of NURBS control points, is interpreted as a Dirichlet boundary condition. The following problem governing the static deformation of $\Omega$ is considered
\begin{equation} \label{eq:elasticityStrong}
	\left\{\begin{aligned}
	{\grad_\bX}\cdot\stress  + \bm{f} & = \0  &&\text{in } \Omega    \\
	\bd(\bX) & = \delta\bd(\bX,\bmu)  &&\text{on } \partial\Omega, 
	\end{aligned} \right. 
\end{equation}
where $\bm{f}$ is an external force defined by the user. The stress tensor $\stress$ is given by
\begin{equation*}
  \stress = \frac{E\nu}{(1+\nu)(1-2\nu)} \tr(\deform)\mat{I} + \frac{E}{1+\nu} \deform ,
\end{equation*}
where $E$ and $\nu$ denote the Young's Modulus and the Poisson's ratio of the elastic medium and the deformation tensor is defined as 
\begin{equation*} 
  \deform = \grad_\bX^S \bd := \frac{1}{2}\left( {\grad_\bX} \bd + ({\grad_\bX} \bd)^T \right) .
\end{equation*}
with $\bd$ being the unknown displacement field.
	
This strategy has been successfully applied in the context of high-order curved mesh generation, see~\cite{Persson2009,HO-Meshing,poya2016unified} for further details. The key aspect here is to write the approximated solution, $\bd_h$, in separated form by using the separated representation of the imposed boundary displacement derived in Section \ref{sbc:boundaryDisplacement}.

The discretisation of the weak formulation associated to the strong form of the solid mechanics problem~\eqref{eq:elasticityStrong} leads to a system of linear equations that can be written as
\begin{equation} \label{eq:solidMechSystemEq}
  \begin{bmatrix}
	\A_{11} & \A_{12} \\
	\A_{21} & \A_{22} \\
  \end{bmatrix}
  \begin{Bmatrix}
	\vect{d} \\
	\vect{\delta d} \\
  \end{Bmatrix}
  =
  \begin{Bmatrix}
	\0 \\
	\0 \\
  \end{Bmatrix},
\end{equation}
where $\vect{d}$ and $\vect{\delta d}$ are vectors containing the nodal values of the approximated displacement $\bd_h$ and the imposed displacement of boundary nodes respectively. 

The matrix $\A_{11}$ in~\eqref{eq:solidMechSystemEq} is symmetric and positive definite. Thus, solving for $\vect{d}$ induces a linear application applied on $\vect{\delta d}$, see \ref{app:solidMechImplementation} for more details. Since $\vect{\delta d}$ has a separable expression, see \eqref{eq:boundayDisplacementSeparated}, the solution of the previous system will induce a separated expression for all nodal values and, consequently, the piecewise interpolation, standard in finite elements, produces a separated representation at any point of the domain, namely
\begin{equation} \label{eq:separatedDisplacement}
 \bd_h(\bX,\bmu) = \sum_{m=1}^N \psi_m^{\bd}(\bmu) \bd_m(\bX) ,
\end{equation}
which is exactly the desired structure of the displacements.

%______________________________________________________________________
\section{Numerical examples} \label{sc:examples}
%______________________________________________________________________

This section presents three numerical examples that show the optimal approximation properties of the proposed PGD approach and its potential for two and three dimensional problems involving geometric parameters as extra coordinates. The examples involve the simulation of Stokes flows using a separable expression that employs the same parametric function for both velocity and pressure. This alternative was shown to be superior to other approaches in\ \cite{diez2017generalized}. For instance,  to satisfy the so-called Ladyzhenskaya-Babu{\v s}ka-Brezzi  (LBB) condition\ \cite{donea2003finite}. Namely, in all the examples $k$ denotes the degree of approximation used for the velocity field and the parametric functions, whereas a degree of approximation $k-1$ is used for the pressure field. This ensures satisfaction of the LBB.

%_________________________________
\subsection{Rotating Couette flow}
%_________________________________

\begin{figure}[!tb]
	\centering
	\includegraphics[width=0.6\textwidth]{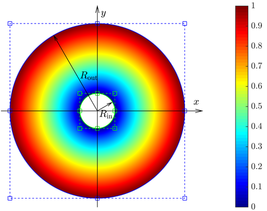}
	\caption{Computational domain for the solution of the rotating Couette flow and magnitude of the velocity of the analytical solution.}
	\label{fig:couetteSetup}
\end{figure}

The first example considers the Couette flow around two infinite coaxial circular cylinders centred at the origin and with radius $\Rin$ and $\Rout$ respectively, with $\Rin < \Rout$, as represented in Figure~\ref{fig:couetteSetup}. The boundary conditions correspond to known angular velocities, $\Vin$ and $\Vout$, at $\Rin$ and $\Rout$, respectively. 
It is worth noting that the pressure must be specified at a point to remove its indeterminacy, as only Dirichlet boundary conditions for the velocity are considered. Here the pressure is imposed at one point of the outer boundary.

The analytical solution for this problem is known~\cite{childs2010rotating}. The azimuthal component of the velocity is given by
\begin{equation} \label{eq:couetteAnalytical}
v_{\phi} = \frac{ \Rout^2 \Vout -\Rin^2 \Vin }{ \Rout^2 -\Rin^2 } r + \frac{ ( \Vin - \Vout) \Rout^2 \Rin^2 }{ \Rout^2 -\Rin^2 } \frac{1}{r}
\end{equation}
where $r = \| \bm{x} \|_2$. The magnitude of the velocity of the analytical solution is also depicted in Figure~\ref{fig:couetteSetup} for $\Vin=0$ and $\Vout = 1$.

The inner radius $\Rin$ is considered an extra parameter within the proposed PGD framework and the objective is to find, in the off-line stage, the generalised velocity and pressure fields for $\Rin \in [1, 2.5]$. It is worth noting that the variation of the inner radius induces the variation of the eight control points (i.e. 16 parameters in two dimensions)  of the NURBS curve describing the inner circle, as represented in Figure~\ref{fig:couetteSetup}. However, as mentioned in Remark~\ref{rk:reparametrisation}, the motion of these control points can be controlled by a single parameter $\mu$ representing the variation of the radius of the inner circle.

\begin{figure}[!tb]
	\centering
	\subfigure[Mesh 1]{\includegraphics[width=0.28\textwidth]{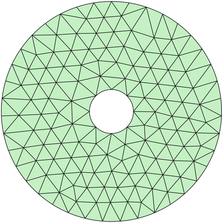}}
	\hspace{0.5cm}
	\subfigure[Mesh 2]{\includegraphics[width=0.28\textwidth]{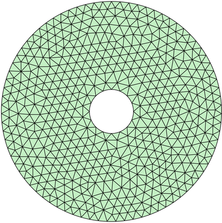}}
	\hspace{0.5cm}
	\subfigure[Mesh 3]{\includegraphics[width=0.28\textwidth]{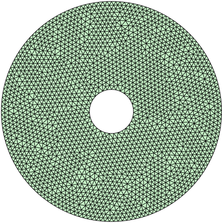}}
	\caption{Three unstructured triangular meshes of the reference domain used for the solution of the Couette flow problem.}
	\label{fig:couetteRefMesh}
\end{figure}

In this example the reference configuration corresponds to $\Rin = 1$ and $\Rout = 5$ and $\mu \in  \I_1 = [0,1.5]$. Three unstructured triangular meshes of the reference domain, with 251, 1,023 and 4,256 elements respectively, are represented in Figure~\ref{fig:couetteRefMesh}. These meshes are generated using the technique proposed in~\cite{NEFEMmeshes} to guarantee that elements without an edge on a curved boundary can be mapped to a reference triangle using an affine mapping.

\begin{figure}[!tb]
	\centering
	\subfigure[$\nDeg=2$]{\includegraphics[width=0.32\textwidth]{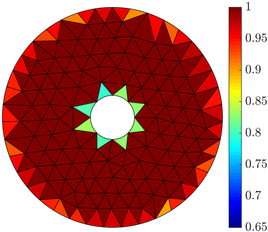}}
	\subfigure[$\nDeg=3$]{\includegraphics[width=0.32\textwidth]{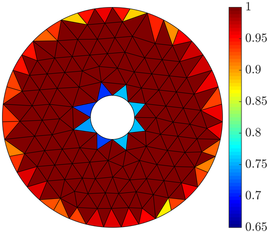}}
	\subfigure[$\nDeg=4$]{\includegraphics[width=0.32\textwidth]{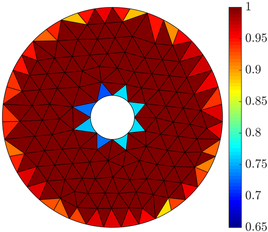}}
	\caption{Quality of the coarsest mesh of the reference domain, shown in Figure~\ref{fig:couetteRefMesh} (a), using different degrees of approximation.}
	\label{fig:couetteRefMeshQuality}
\end{figure}

The quality of the coarsest mesh of the reference domain, measured using the scaled Jacobian~\cite{HO-Meshing}, is represented in Figure~\ref{fig:couetteRefMeshQuality}. 
The minimum quality observed is, as expected near the inner boundary where the elements show the maximum distortion. 

\begin{figure}[!tb]
	\centering
	\subfigure[$\mu=0.5$]{\includegraphics[width=0.32\textwidth]{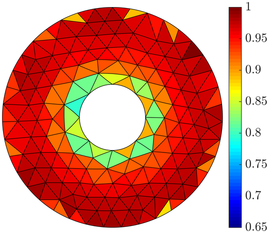}}
	\subfigure[$\mu=1$]  {\includegraphics[width=0.32\textwidth]{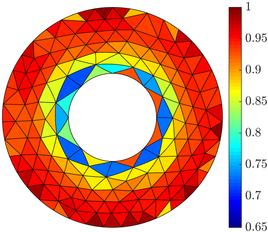}}
	\subfigure[$\mu=1.5$]{\includegraphics[width=0.32\textwidth]{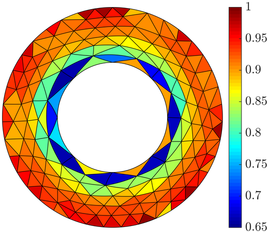}}
	\caption{Quality of the coarsest mesh of the deformed domain with $\nDeg=4$ for different values of the parameter $\mu$.}
	\label{fig:couetteDefMeshQualityP4}
\end{figure}

Figure~\ref{fig:couetteDefMeshQualityP4} shows the deformed configuration for three values of the parameter $\mu$, namely $\mu=0.5$, $\mu=1$ and $\mu=1.5$. These plots also represent the quality of the deformed mesh obtained by using the elastic analogy described in Section~\ref{sbc:domainDisplacement}.
In all cases, the minimum quality is higher than 0.65.

\begin{figure}[!tb]
	\centering
	\subfigure[$\nDeg=2$]{\includegraphics[width=0.32\textwidth]{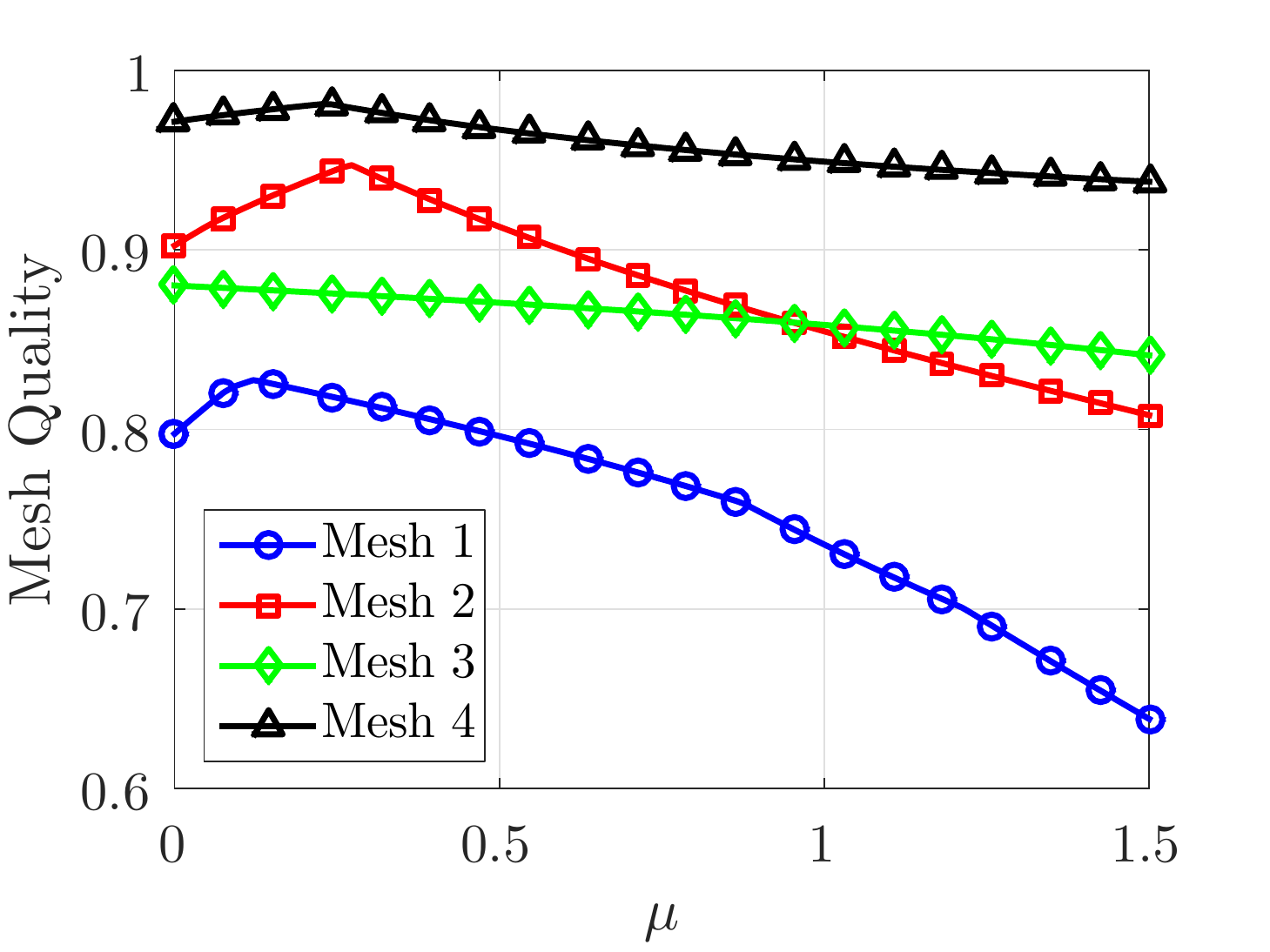}}
	\subfigure[$\nDeg=3$]{\includegraphics[width=0.32\textwidth]{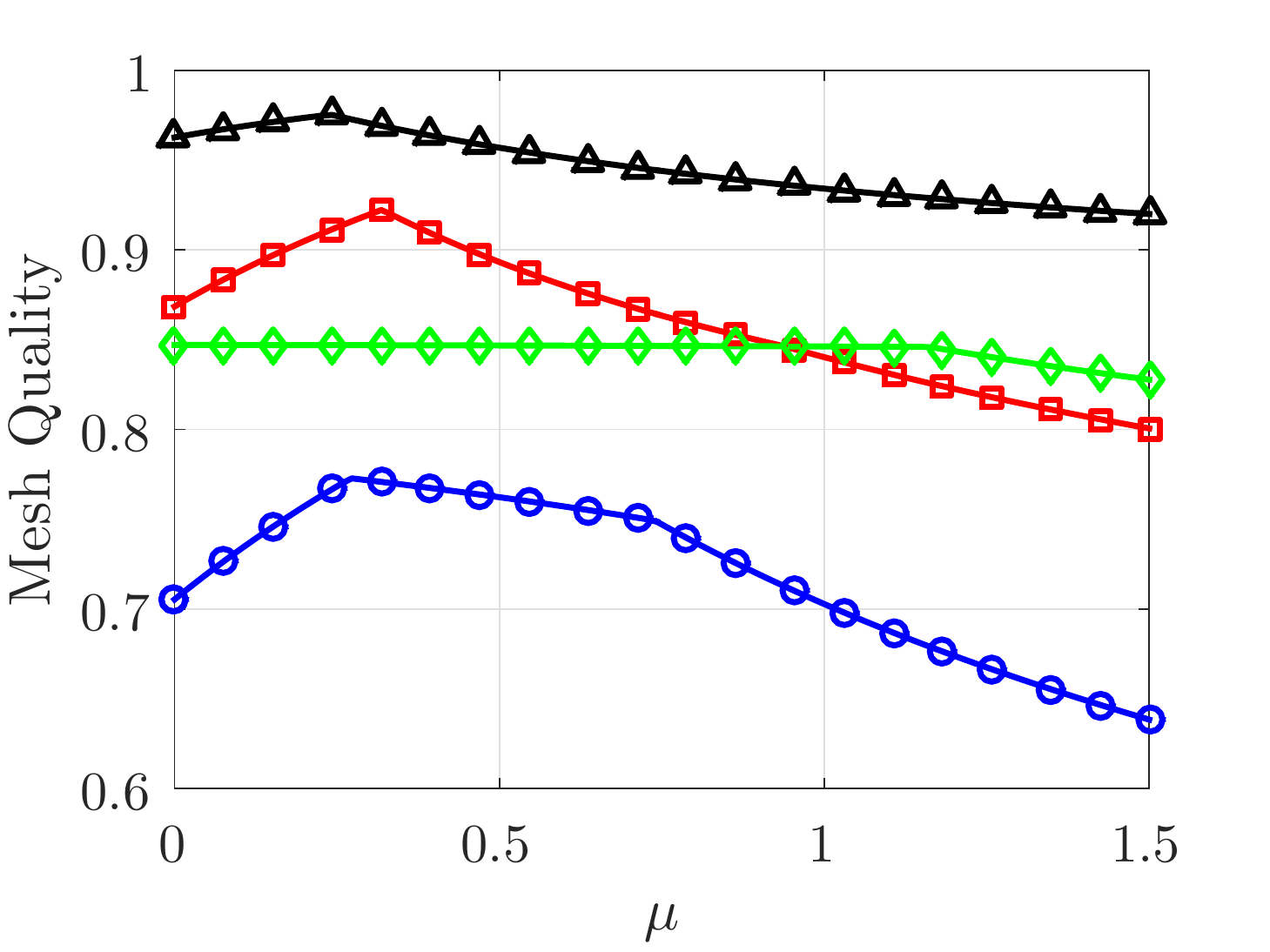}}
	\subfigure[$\nDeg=4$]{\includegraphics[width=0.32\textwidth]{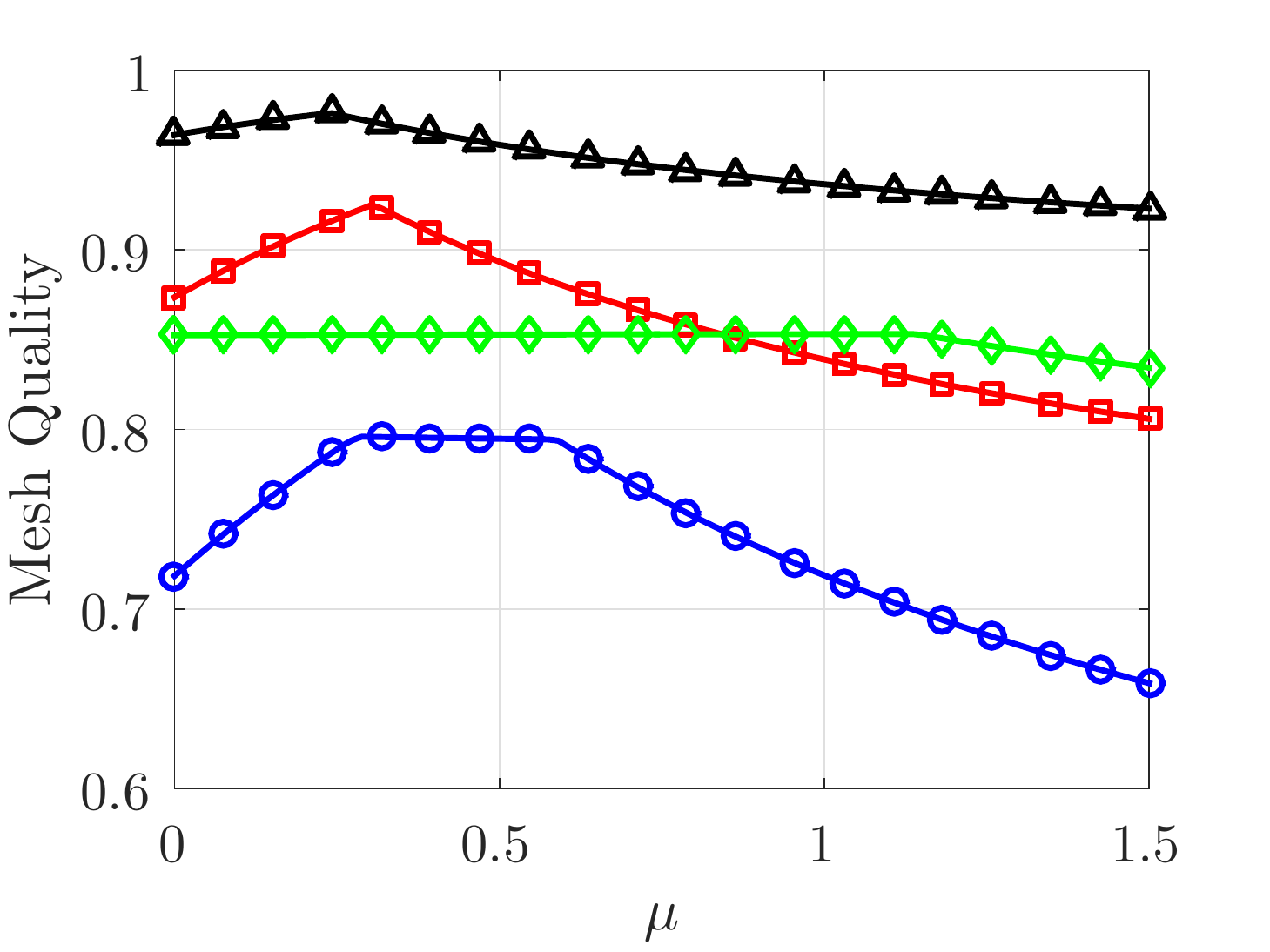}}
	\caption{Evolution of the quality of the deformed meshes corresponding to the reference mesh of Figure~\ref{fig:couetteRefMesh} (a) as a function of the parameter $\mu$ using different meshes and degrees of approximation.}
	\label{fig:couetteDefMeshQuality}
\end{figure}

To further illustrate the robustness of the mesh deformation technique employed within the proposed PGD framework, Figure~\ref{fig:couetteDefMeshQuality} shows the evolution of the quality of the deformed meshes as a function of the parameter $\mu$ using four different meshes and three different degrees of approximation.
As expected, the worst case scenario corresponds to the maximum deformation induced by a parameter $\mu=1.5$, but in all cases the minimum quality is always higher than 0.6. Furthermore, it is worth noting that the influence of the parameter $\mu$ on the quality of the deformed meshes is less important for finer meshes.

A crucial aspect of the proposed PGD strategy is the separation of the matrix $\mat{H}_{\bmu}$ defined in Equation~\eqref{eq:H}. As discussed previously, $\mat{H}_{\bmu}$ does not generally admit an exact separable expression and therefore a separable approximation is computed here via the higher-order PGD-projection~\cite{DM-MZH:15}. 

\begin{figure}[!tb]
	\centering
	\subfigure[$m=1$]{\includegraphics[width=0.24\textwidth]{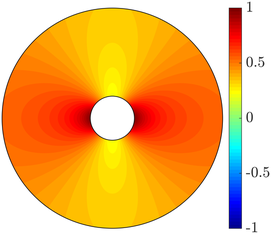}}
	\subfigure[$m=2$]{\includegraphics[width=0.24\textwidth]{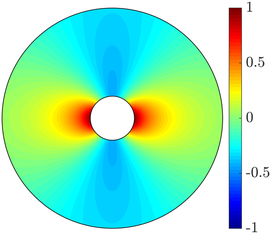}}
	\subfigure[$m=3$]{\includegraphics[width=0.24\textwidth]{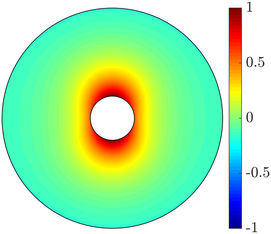}}
	\subfigure[$m=4$]{\includegraphics[width=0.24\textwidth]{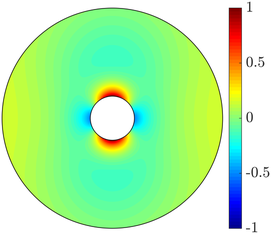}}
	\subfigure[$m=5$]{\includegraphics[width=0.24\textwidth]{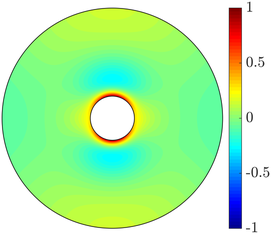}}
	\subfigure[$m=6$]{\includegraphics[width=0.24\textwidth]{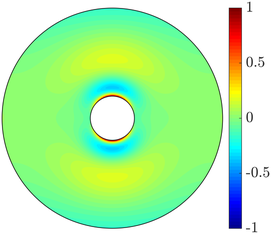}}
	\subfigure[$m=7$]{\includegraphics[width=0.24\textwidth]{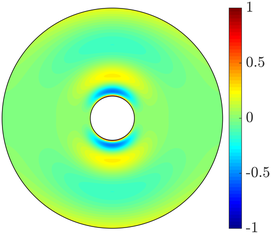}}
	\subfigure[$m=8$]{\includegraphics[width=0.24\textwidth]{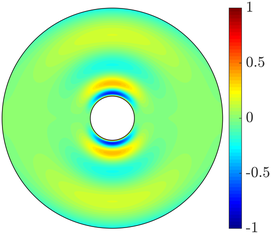}}
	\caption{First eight normalised spatial modes of $[\mat{H}_{\bmu}]_{11}$ on the mesh of Figure~\ref{fig:couetteRefMesh} (c) with 800 elements in the parametric dimension and with $\nDeg=4$.}
	\label{fig:couetteMesh_H3P4_H11_Mode}
\end{figure}

\begin{figure}[!tb]
	\centering
	\subfigure[$m=1$]{\includegraphics[width=0.24\textwidth]{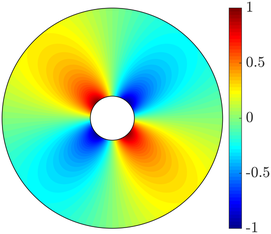}}
	\subfigure[$m=2$]{\includegraphics[width=0.24\textwidth]{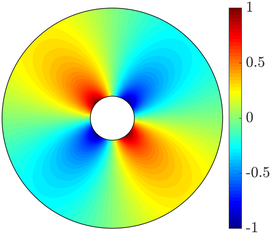}}
	\subfigure[$m=3$]{\includegraphics[width=0.24\textwidth]{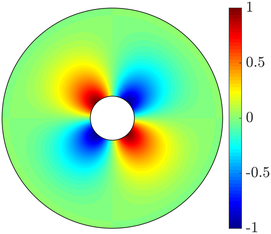}}
	\subfigure[$m=4$]{\includegraphics[width=0.24\textwidth]{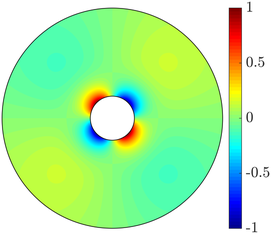}}
	\subfigure[$m=5$]{\includegraphics[width=0.24\textwidth]{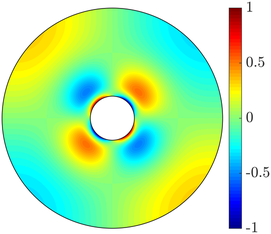}}
	\subfigure[$m=6$]{\includegraphics[width=0.24\textwidth]{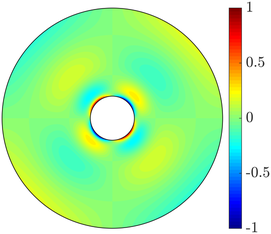}}
	\subfigure[$m=7$]{\includegraphics[width=0.24\textwidth]{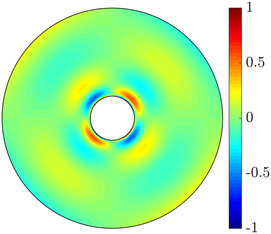}}
	\subfigure[$m=8$]{\includegraphics[width=0.24\textwidth]{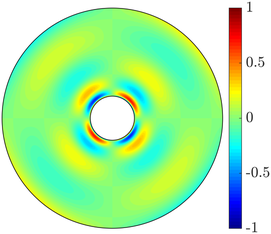}}
	\caption{First eight normalised spatial modes of $[\mat{H}_{\bmu}]_{12}$ on the mesh of Figure~\ref{fig:couetteRefMesh} (c) with 800 elements in the parametric dimension and with $\nDeg=4$.}
	\label{fig:couetteMesh_H3P4_H12_Mode}
\end{figure}

Figure~\ref{fig:couetteMesh_H3P4_H11_Mode} shows the first eight normalised spatial modes of the component $[\mat{H}_{\bmu}]_{11}$.
The results suggest that the first two modes capture the main global features of $[\mat{H}_{\bmu}]_{11}$, whereas the rest of modes capture the local variations near the inner circle. The results for the second diagonal component, $[\mat{H}_{\bmu}]_{22}$, not displayed for brevity, show the same behaviour but with the expected rotation of 90 degrees due to the symmetry of the domain and the displacement field. Similarly, Figure~\ref{fig:couetteMesh_H3P4_H12_Mode} shows the first eight normalised spatial modes of the component $[\mat{H}_{\bmu}]_{12}$.
The results show, again how the first modes capture the global behaviour of the component $[\mat{H}_{\bmu}]_{12}$, whereas the last modes show relevant spatial variations in the close vicinity of the inner circle. 

\begin{figure}[!tb]
	\centering
	\includegraphics[width=0.55\textwidth]{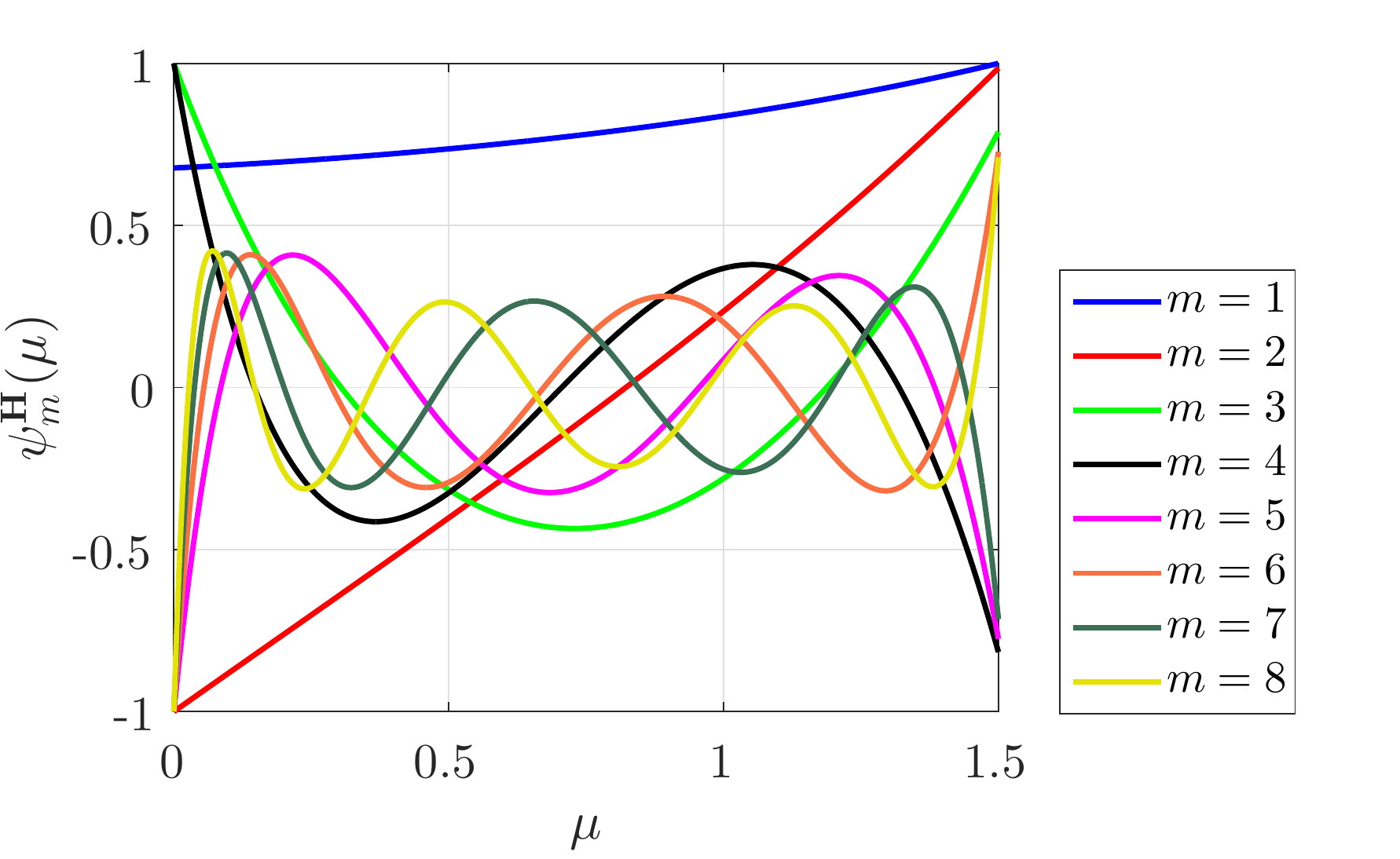}
	\caption{First eight parametric modes of $\mat{H}_{\bmu}$ on the mesh of Figure~\ref{fig:couetteRefMesh} (c) with 800 elements in the parametric dimension and with $\nDeg=4$.}
	\label{fig:couetteMesh_H3P4_HXX_ParamModes}
\end{figure}

The first eight normalised parametric modes of $\mat{H}_{\bmu}$ are represented in Figure~\ref{fig:couetteMesh_H3P4_HXX_ParamModes}. 
It is worth recalling that the same parametric modes are associated to all the components of the matrix $\mat{H}_{\bmu}$.

\begin{figure}[!tb]
	\centering
	\subfigure[$\nDeg=2$]{\includegraphics[width=0.32\textwidth]{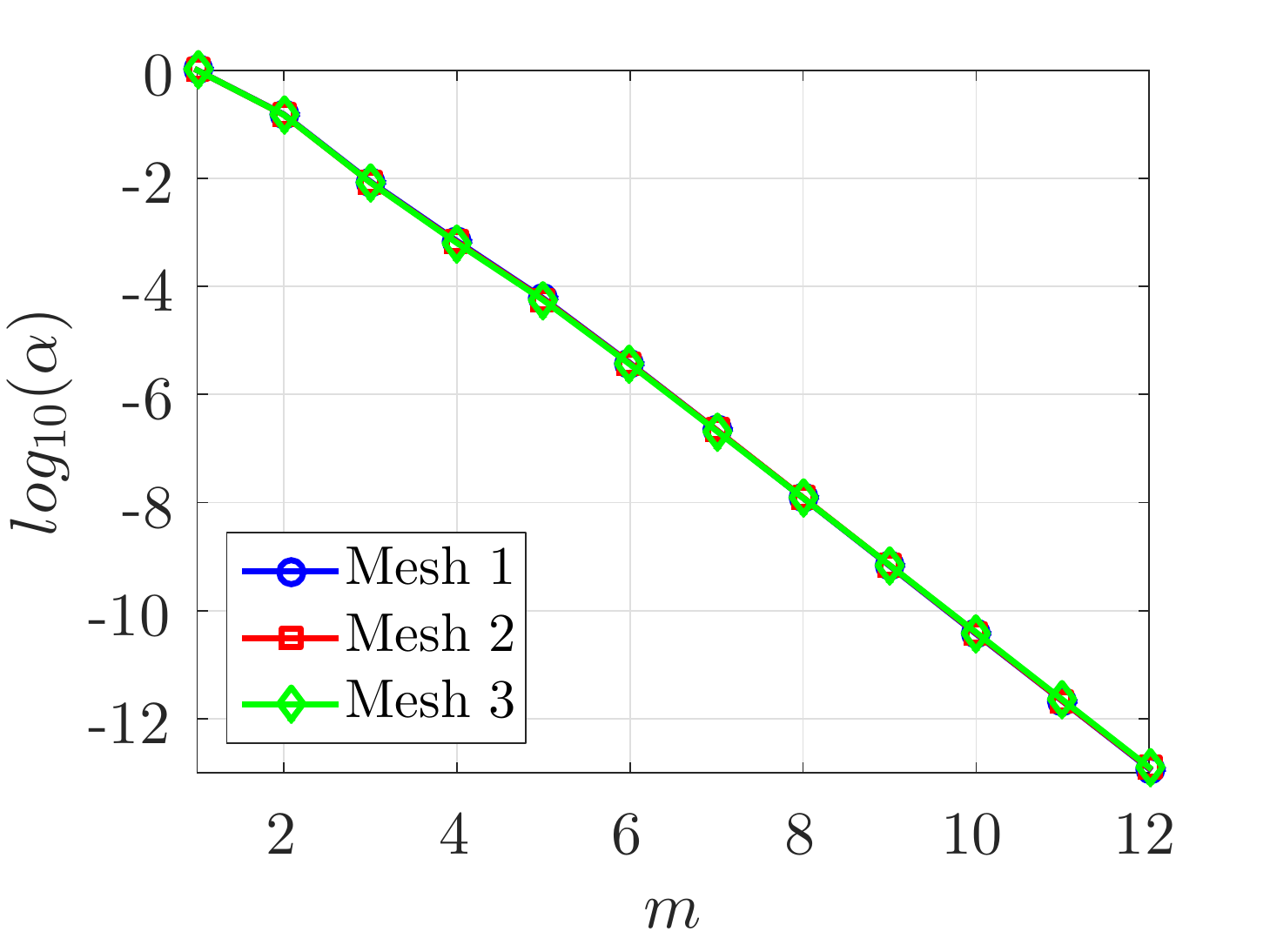}}
	\subfigure[$\nDeg=3$]{\includegraphics[width=0.32\textwidth]{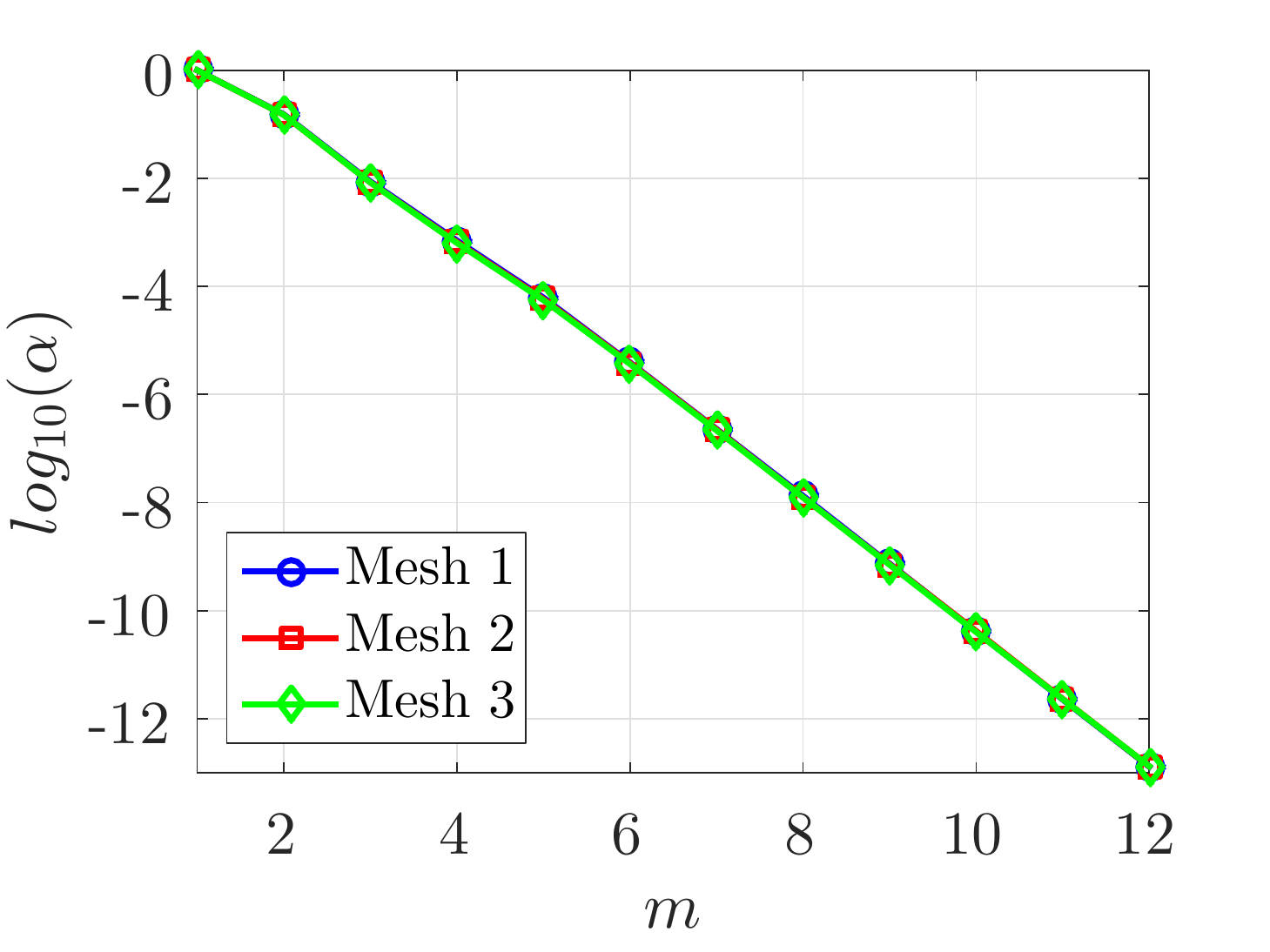}}
	\subfigure[$\nDeg=4$]{\includegraphics[width=0.32\textwidth]{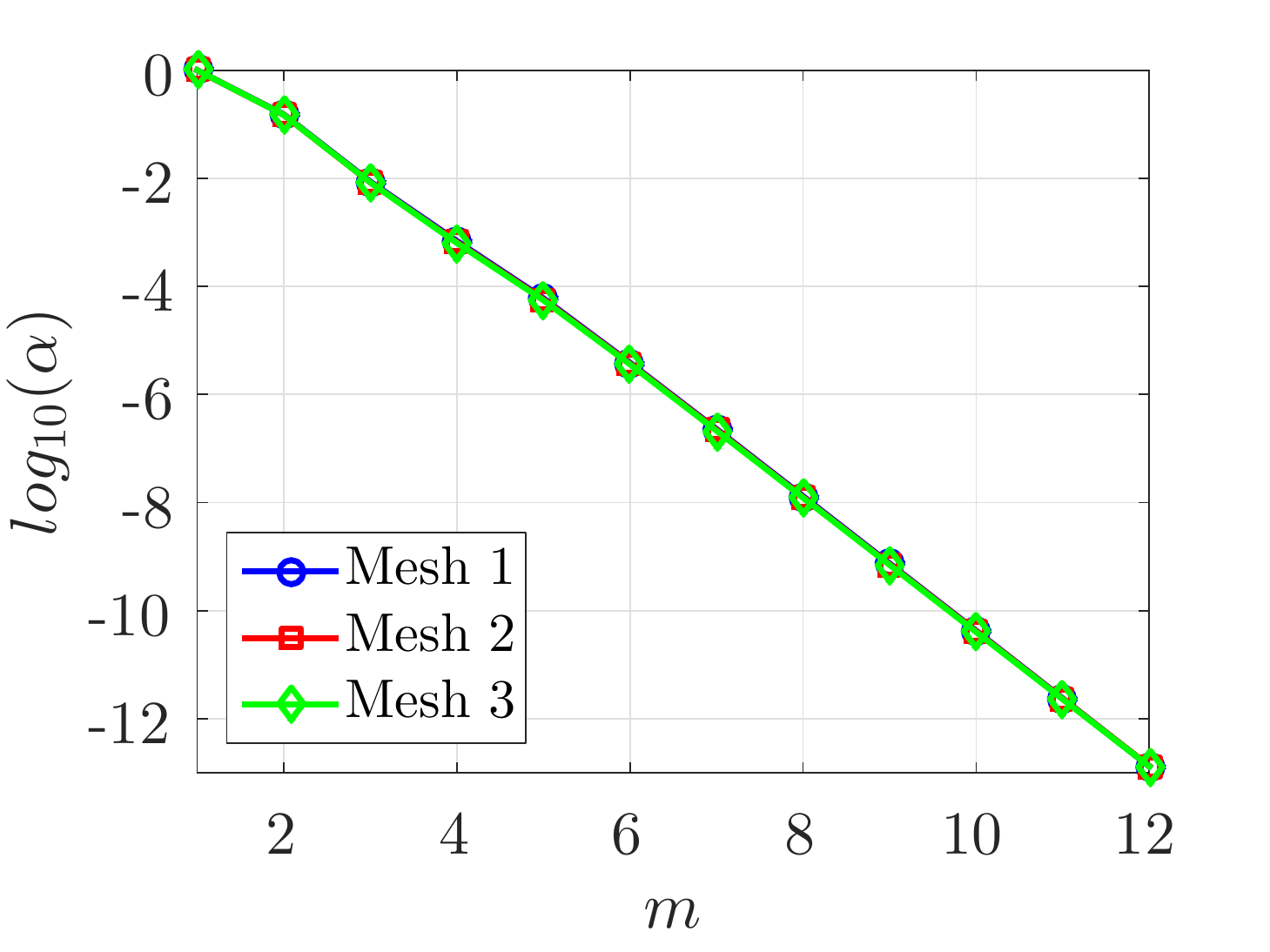}}
	\caption{Evolution of the amplitude of the spatial modes $\alpha_m$ of the matrix $\mat{H}_{\bmu}$ using different meshes and degrees of approximation.}
	\label{fig:couette_AmplitudeH11}
\end{figure}

The amplitude, $\alpha_m$, corresponding to the mode $m$ of the separation of $\mat{H}_{\bmu}$, is computed as the product of the Euclidean norms of the spatial and parametric functions. Figure~\ref{fig:couette_AmplitudeH11} shows the amplitudes of $[\mat{H}_{\bmu}]_{11}$ using three different meshes and three different degrees of approximation.
The results show that, for this example, the number of terms required to obtain a separable approximation of the matrix $\mat{H}_{\bmu}$ using the higher-order PGD-projection is completely independent on the spatial discretisation. In all cases, 12 modes provide a decrease in the amplitude $\alpha_m$ of exactly 13 orders of magnitude and the amplitude is the same in all meshes and for all degrees of approximation. 

\begin{figure}[!tb]
	\centering
	\subfigure[$m=1$]{\includegraphics[width=0.24\textwidth]{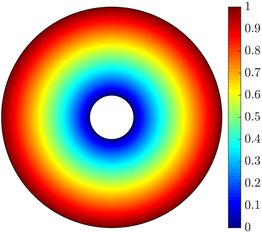}}
	\subfigure[$m=2$]{\includegraphics[width=0.24\textwidth]{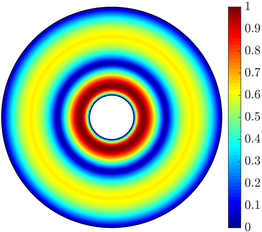}}
	\subfigure[$m=3$]{\includegraphics[width=0.24\textwidth]{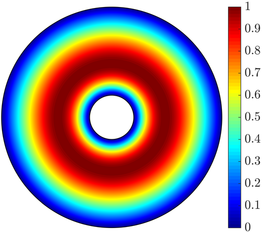}}
	\subfigure[$m=4$]{\includegraphics[width=0.24\textwidth]{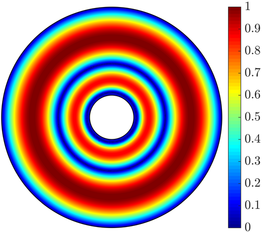}}
	\subfigure[$m=5$]{\includegraphics[width=0.24\textwidth]{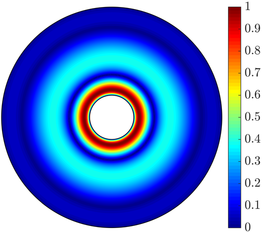}}
	\subfigure[$m=6$]{\includegraphics[width=0.24\textwidth]{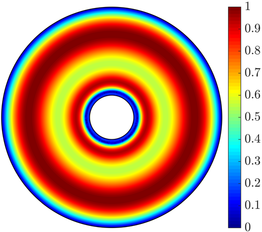}}
	\subfigure[$m=7$]{\includegraphics[width=0.24\textwidth]{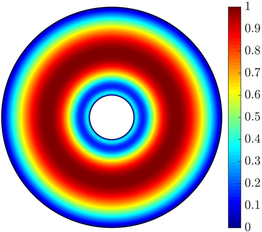}}
	\subfigure[$m=8$]{\includegraphics[width=0.24\textwidth]{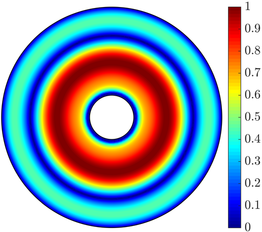}}
	\caption{First eight normalised spatial modes of the magnitude of the velocity computed on the mesh of Figure~\ref{fig:couetteRefMesh} (c),  800 elements in the parametric dimension and a degree of approximation $\nDeg=4$.}
	\label{fig:couetteModesV}
\end{figure}

Using the separation of the matrix $\mat{H}_{\bmu}$, the rotating Couette flow problem is solved to obtain the generalised solution of the Stokes problem. The first eight normalised spatial modes of the magnitude of the velocity field are shown in Figure~\ref{fig:couetteModesV}. The simulation was performed using the mesh of Figure~\ref{fig:couetteRefMesh} (c), 800 equally-spaced elements in $\I_1 = [0,1.5]$ and a degree of approximation $\nDeg=4$. 
The solution of the coaxial Couette flow problem requires imposing Dirichlet boundary conditions for the velocity at the inner and outer boundaries of the spatial domain. This is implemented within the PGD framework by adding one initial mode that fulfils the non-homogeneous Dirichlet boundary conditions. Its spatial part is computed as the solution of the Stokes problem at the reference mesh and its parametric component is taken as a constant function equal to one. The next PGD modes are computed by imposing homogeneous Dirichlet boundary conditions, ensuring that the separated PGD solution satisfies the required Dirichlet boundary conditions.

\begin{figure}[!tb]
	\centering
	\includegraphics[width=0.55\textwidth]{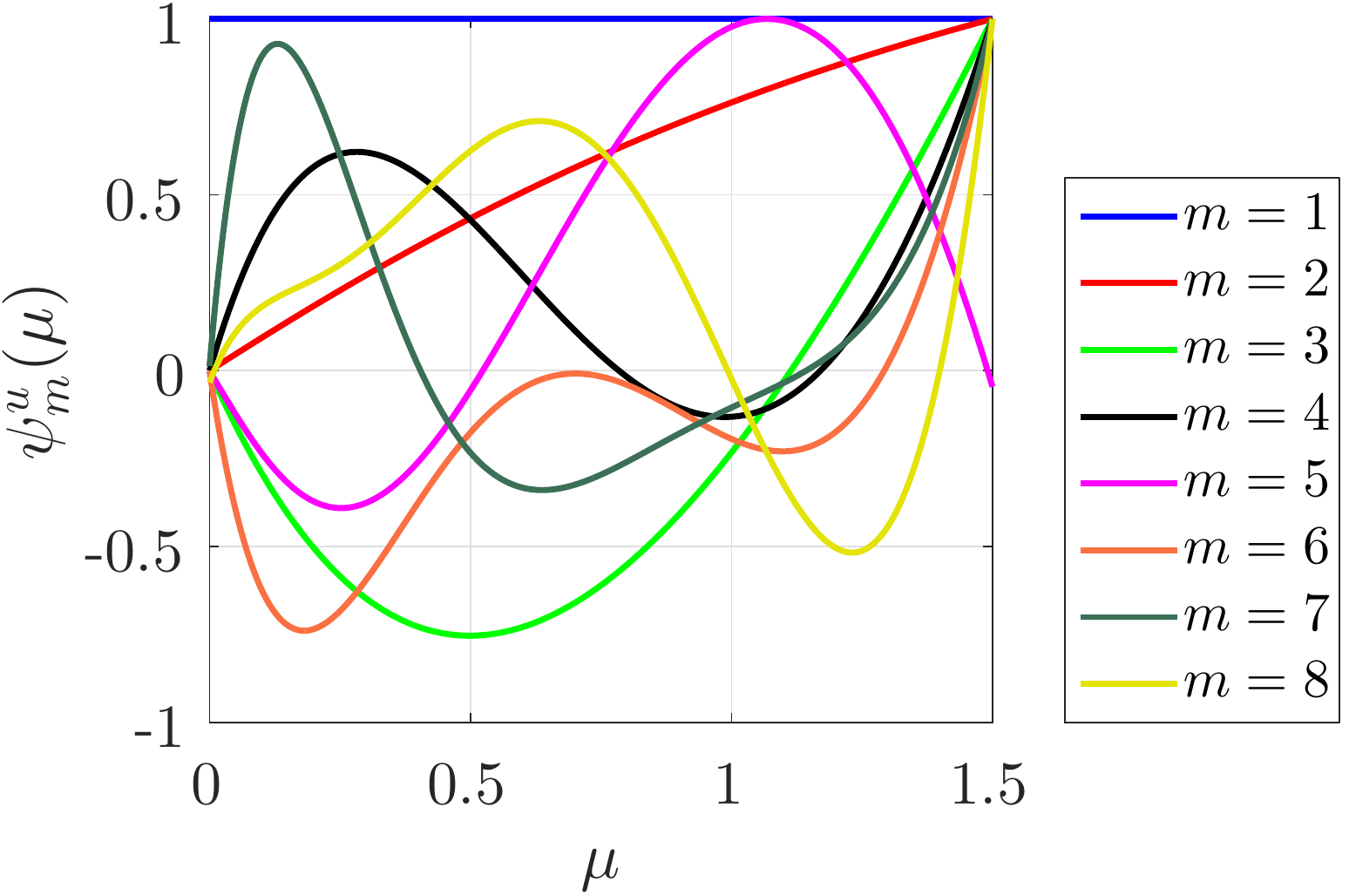}
	\caption{First eight parametric modes of the rotating Couette flow on the mesh of Figure~\ref{fig:couetteRefMesh} (c) with 800 elements in the parametric dimension and with $\nDeg=4$.}
	\label{fig:couette_ParamModes}
\end{figure}

The first eight normalised parametric modes associated to the spatial modes of Figure~\ref{fig:couetteModesV} are represented in Figure~\ref{fig:couette_ParamModes}. 
It is worth recalling that the same parametric mode is associated to all the components of the velocity and the pressure fields. The results reveal that the parametric modes of the velocity field show a similar behaviour compared to the parametric modes of the separation of $\mat{H}_{\bmu}$, shown in Figure~\ref{fig:couetteMesh_H3P4_HXX_ParamModes}. The first two modes are smooth whereas the next modes, that contribute less to the global solution, show a more oscillatory character.

\begin{figure}[!tb]
	\centering
	\subfigure[$\mu=0.5$, $N=1$]{\includegraphics[width=0.32\textwidth]{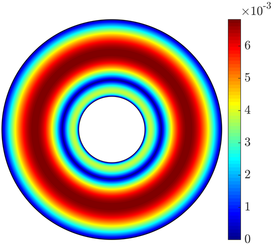}}
	\subfigure[$\mu=0.5$, $N=2$]{\includegraphics[width=0.32\textwidth]{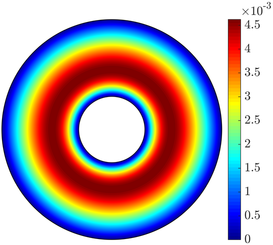}}
	\subfigure[$\mu=0.5$, $N=4$]{\includegraphics[width=0.32\textwidth]{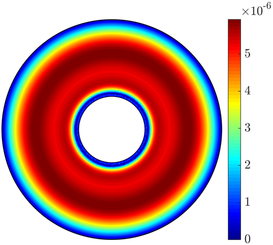}}
	\subfigure[$\mu=1$, $N=1$]{\includegraphics[width=0.32\textwidth]{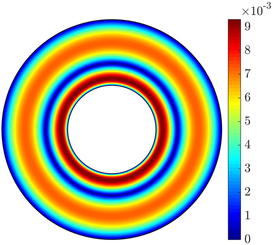}}
	\subfigure[$\mu=1$, $N=2$]{\includegraphics[width=0.32\textwidth]{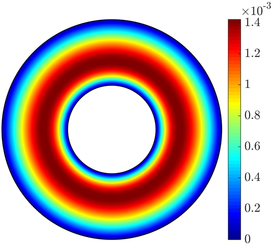}}
	\subfigure[$\mu=1$, $N=4$]{\includegraphics[width=0.32\textwidth]{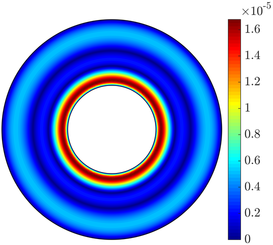}}
	\subfigure[$\mu=1.5$, $N=1$]{\includegraphics[width=0.32\textwidth]{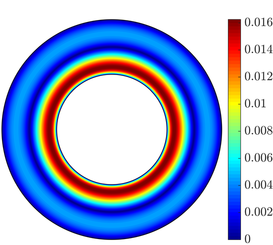}}
	\subfigure[$\mu=1.5$, $N=2$]{\includegraphics[width=0.32\textwidth]{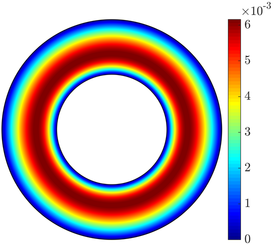}}
	\subfigure[$\mu=1.5$, $N=4$]{\includegraphics[width=0.32\textwidth]{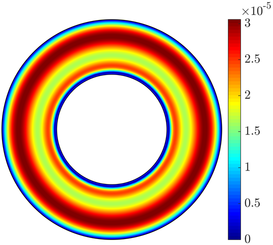}}
	\caption{Absolute value of the error in the magnitude of the velocity computed with the proposed PGD approach using $N$ modes for different values of the parameter $\mu$. The computation has been performed on the mesh of Figure~\ref{fig:couetteRefMesh} (c), with 800 elements in the parametric dimension and a degree of approximation $\nDeg=4$.}
	\label{fig:couette_ErrVelo_H3P4}
\end{figure}

To illustrate the gain in accuracy as the number of modes increases, Figure~\ref{fig:couette_ErrVelo_H3P4} shows the absolute value of the error in the magnitude of the velocity computed with the proposed PGD approach using $N$ modes for different values of the parameter $\mu$. 
The computation has been performed on the mesh of Figure~\ref{fig:couetteRefMesh} (c), with 800 elements in the parametric dimension and a degree of approximation $\nDeg=4$. In all cases, the results illustrate that the computation with four modes is almost three orders of magnitude more accurate than using a single mode. 

\begin{figure}[!tb]
	\centering
	\subfigure[$\nDeg=2$]{\includegraphics[width=0.32\textwidth]{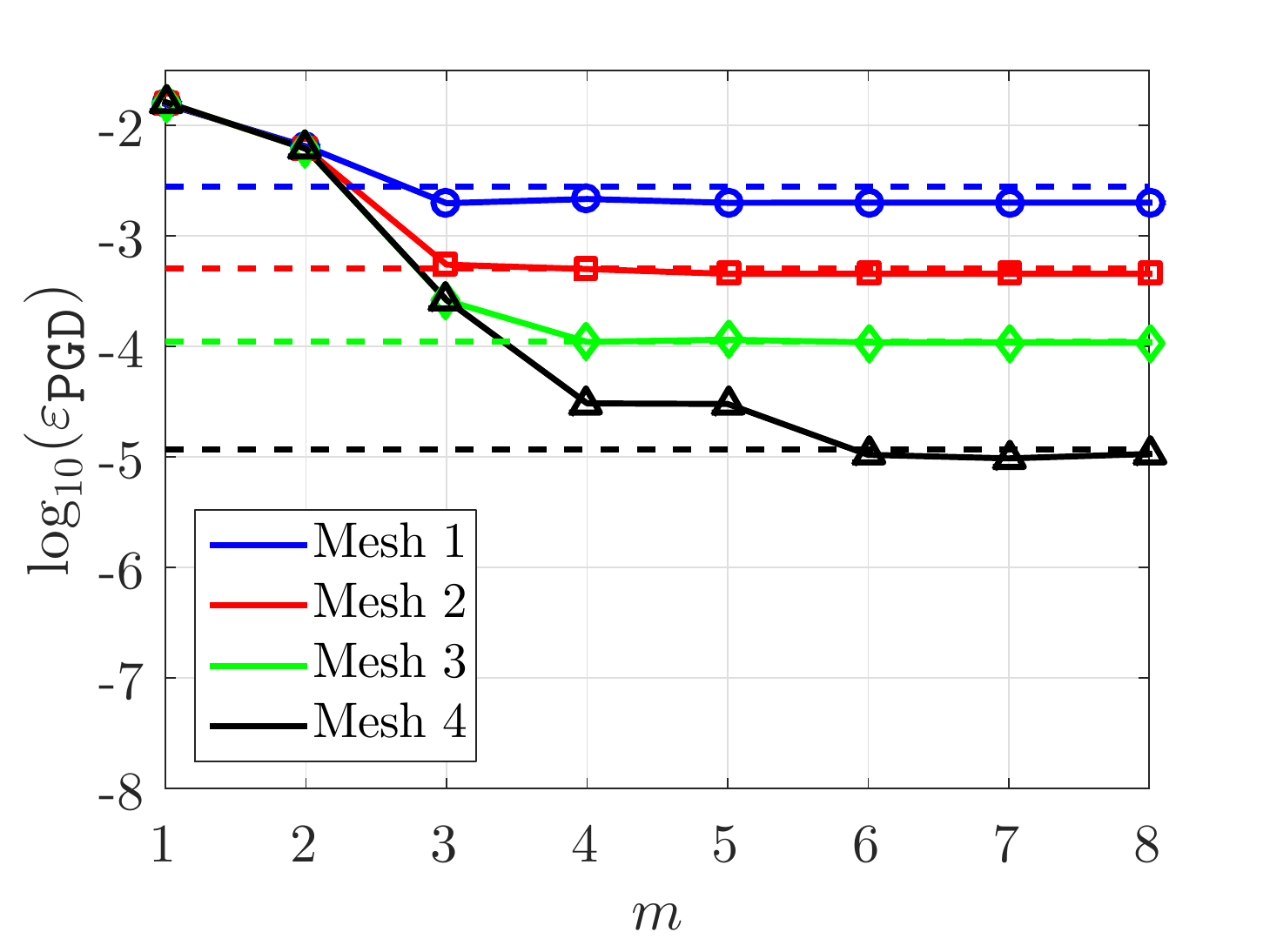}}
	\subfigure[$\nDeg=3$]{\includegraphics[width=0.32\textwidth]{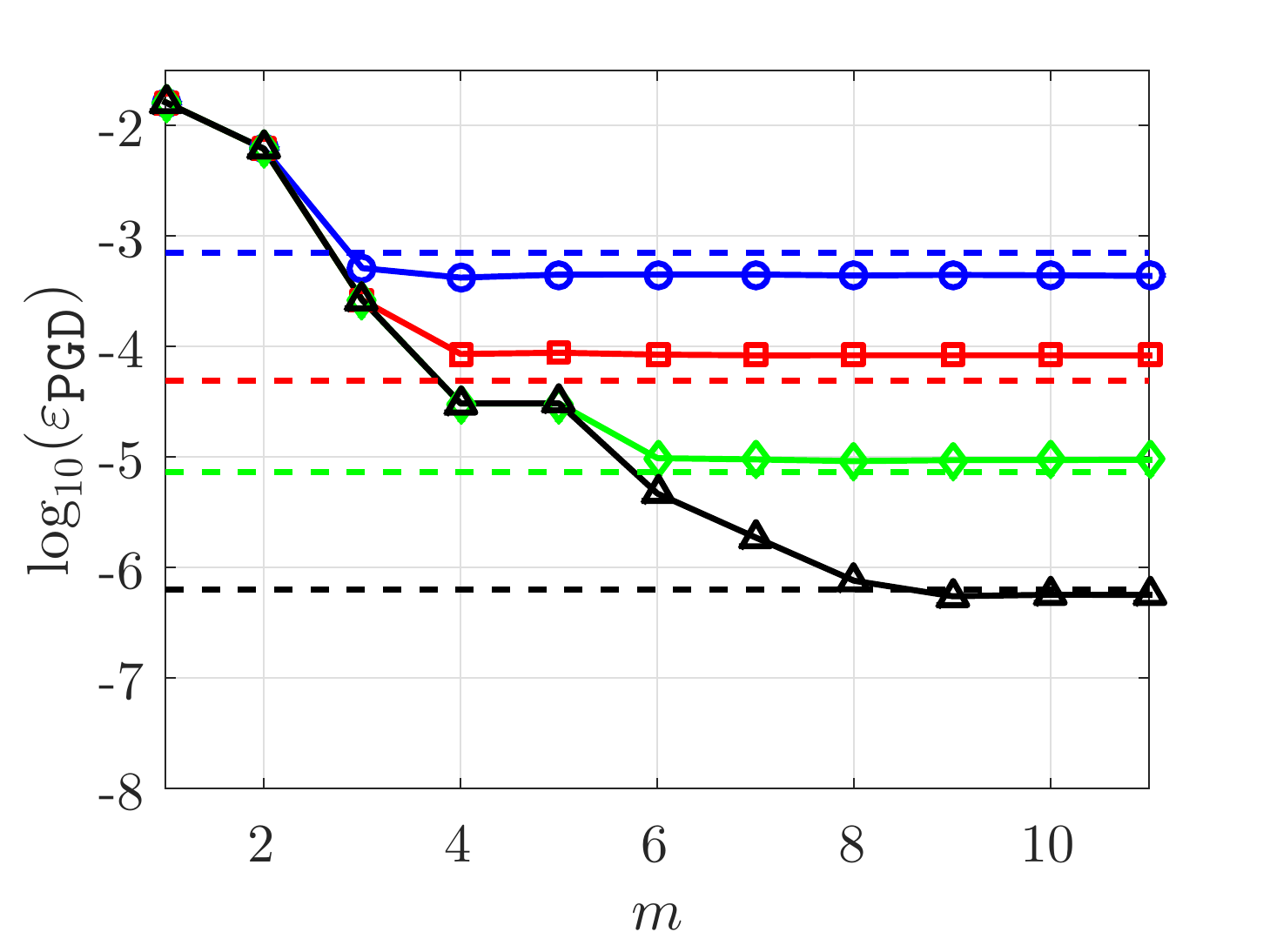}}
	\subfigure[$\nDeg=4$]{\includegraphics[width=0.32\textwidth]{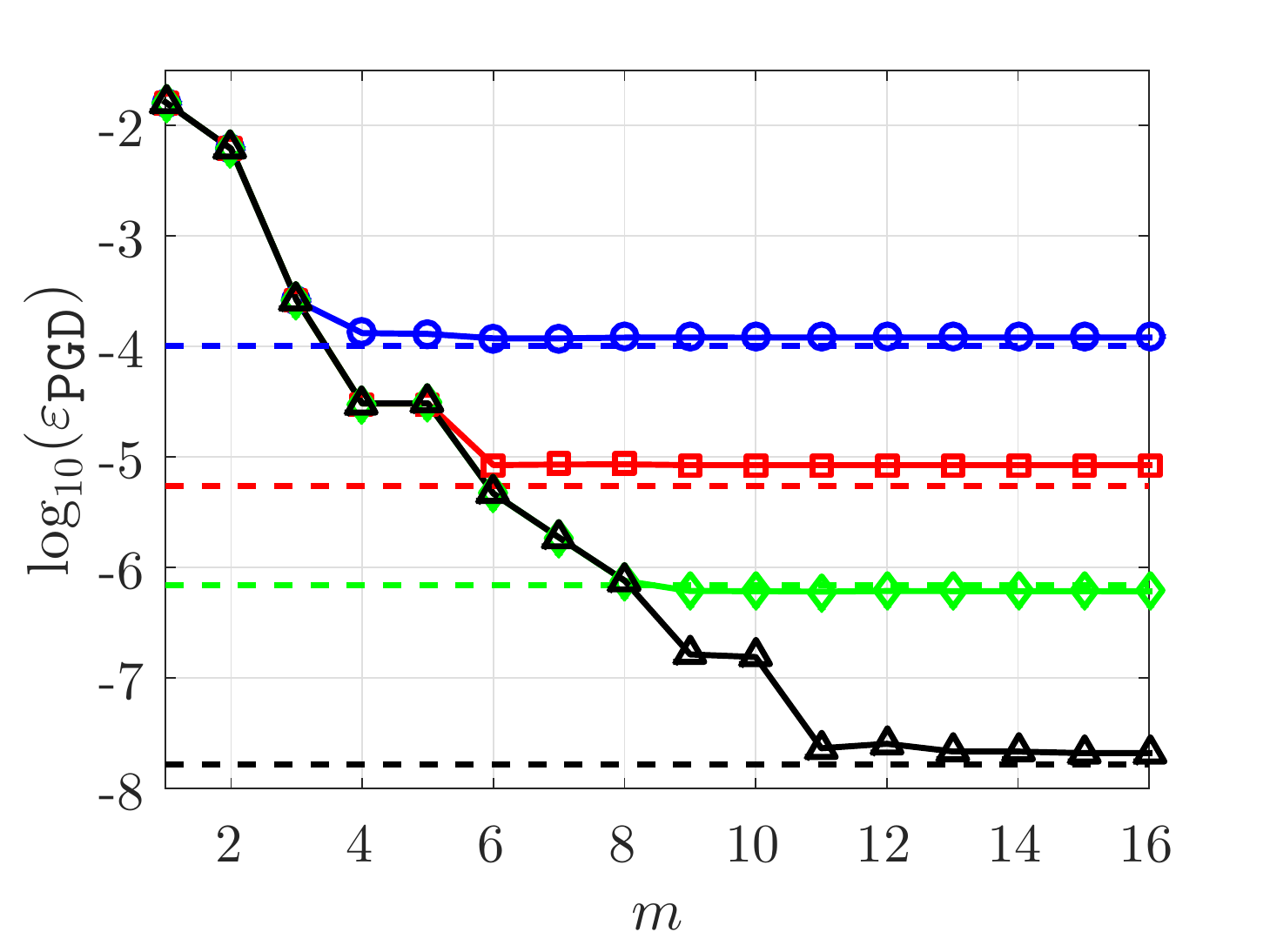}}
	\caption{Error of the PGD solution as a function of the number of PGD modes using three different meshes for the solution of the Couette flow problem. The discontinuous line represents the error using standard FEM.}
	\label{fig:couetteConvModes}
\end{figure}

To further analyse the accuracy of the proposed PGD approach when the number of modes is increased, the relative error in the $\eltwo(\Omega \times \I_1)$ norm, defined as
\begin{equation} \label{eq:PGDerrorAnalytical}
 \varepsilon_{\texttt{PGD}} 
 = \left( \frac{ \displaystyle \int_{\I_1} \int_{\Omega} ( \bm{u}_{\texttt{PGD}} - \bm{u} ) \cdot ( \bm{u}_{\texttt{PGD}} - \bm{u} ) d\Omega \,  d\mu }
                   { \displaystyle \int_{\I_1} \int_{\Omega}  \bm{u} \cdot \bm{u}\, d\Omega \, d\mu} \right)^{1/2} ,
\end{equation}
is studied. Figure~\ref{fig:couetteConvModes} depicts the evolution of $\varepsilon_{\texttt{PGD}}$ as a function of the number of PGD modes for three different degrees of approximation and four different meshes.
The discontinuous lines in Figure~\ref{fig:couetteConvModes} show the relative error of the standard FEM solution measured in the $\eltwo(\Omega \times \I_1)$ norm. The evaluation of the error of the FEM solution in the $\eltwo(\Omega \times \I_1)$ norm requires the solution of as many FEM problems as the number of elements in the parametric dimension used by the PGD approach, multiplied by the number of integration points within each element. 

The results clearly illustrate the increased level of accuracy as the number of modes is increased. In addition, the increased accuracy induced by the use of finer meshes or higher degrees of approximation can be observed. It worth noting that the PGD solution achieves its maximum accuracy (i.e. the accuracy of a standard FEM solution) using a low number of modes. For coarse meshes, the maximum accuracy of the PGD solution is achieved with three or four modes. For finer meshes, the PGD approach provides accurate results with only three or four modes and requires between six and 11 modes to reach the same level of error as the standard FEM solution. 

\begin{figure}[!tb]
	\centering
	\subfigure[$\nDeg=2$]{\includegraphics[width=0.32\textwidth]{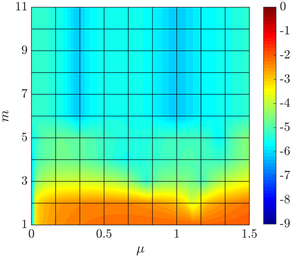}}
	\subfigure[$\nDeg=3$]{\includegraphics[width=0.32\textwidth]{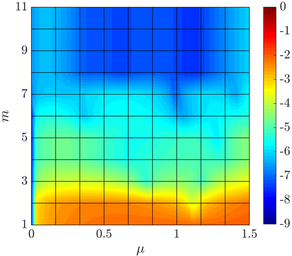}}
	\subfigure[$\nDeg=4$]{\includegraphics[width=0.32\textwidth]{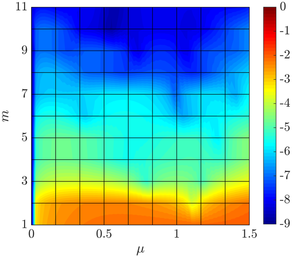}}
	\caption{Difference between the PGD solution and the FEM solution, in logarithmic scale, as a function of the parameter $\mu$ and the number of PGD modes using different degrees of approximation.}
	\label{fig:couetteErrMapParamModes}
\end{figure}

Next, the accuracy of the PGD solution for different values of the parameter $\mu$ is studied. Figure~\ref{fig:couetteErrMapParamModes} shows the difference between the PGD solution and the FEM solution, in logarithmic scale, as a function of the parameter $\mu$ and the number of PGD modes using different degrees of approximation.
The results show that the accuracy of the PGD solution depends only weakly on the value of the parameter selected. For values of $\mu=0$, corresponding to the undeformed configuration, the error is lower for a moderate number of PGD modes. This is expected because for $\mu=0$ the displacement function $\bd$ is null and the Jacobian $\Jaco$ is equal to the identity matrix. But for a large enough number of modes, the difference between the PGD solution and the FEM solution shows very little dependence on the value of $\mu$, illustrating the robustness of the proposed PGD approach.

\begin{figure}[!tb]
	\centering
	\includegraphics[width=0.5\textwidth]{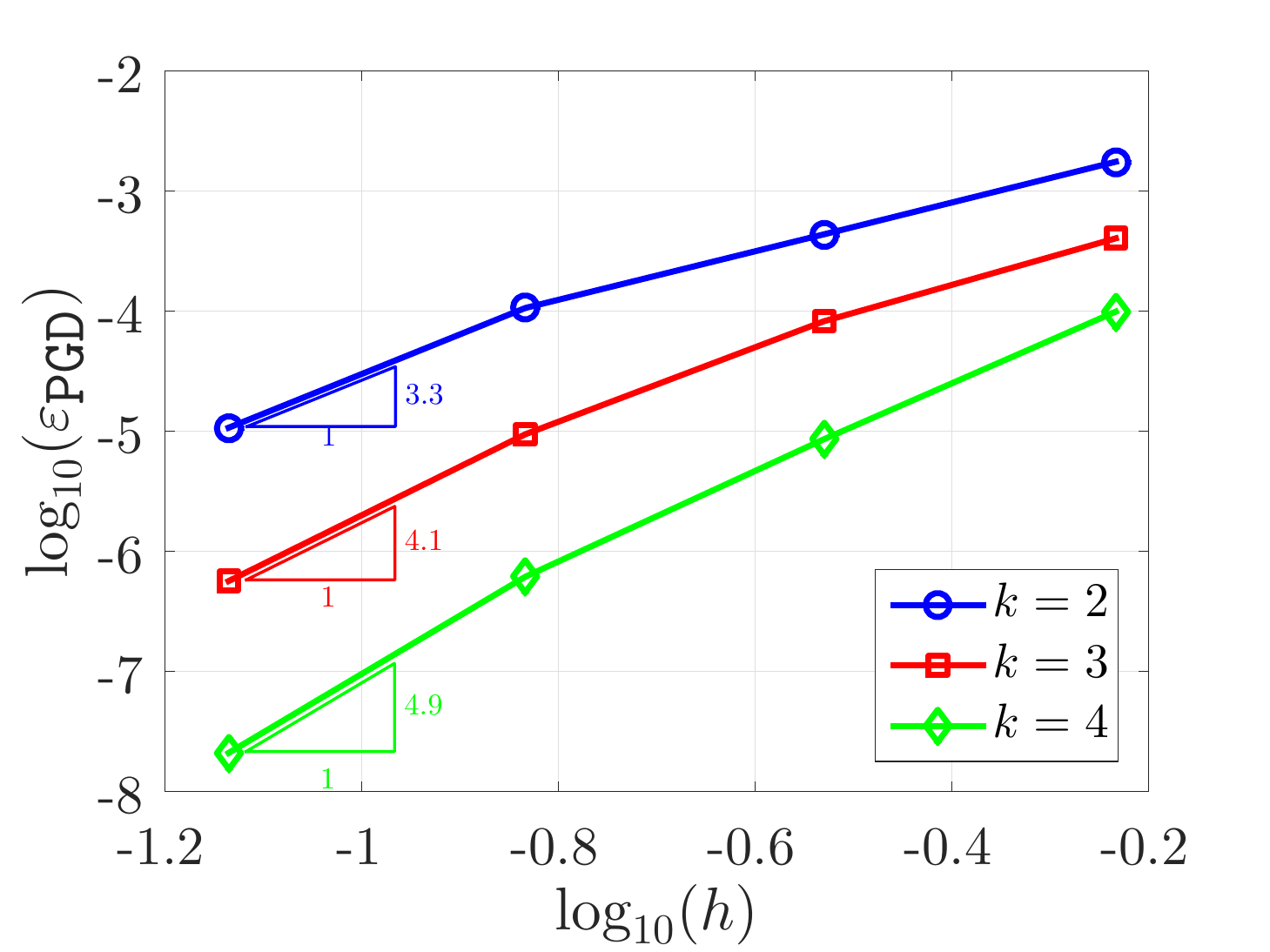}
	\caption{Relative error between the PGD and analytical solutions in the  $\eltwo(\Omega \times \I_1)$ norm as a function of the characteristic element size $h$ for different degrees of approximation $\nDeg$.}
	\label{fig:couetteMesh_hConvPGD}
\end{figure}

The last numerical study for the rotation Couette flow considers a mesh convergence analysis for the proposed PGD approach. The PGD solution is computed on the three meshes shown in Figure~\ref{fig:couetteRefMesh} using three different orders of approximation. Figure~\ref{fig:couetteMesh_hConvPGD} shows the evolution of the relative error in the $\eltwo(\Omega \times \I_1)$ norm, $\varepsilon_{\texttt{PGD}}$, as a function of the characteristic element size $h$.
For each simulation the minimum number of PGD modes required to achieve the maximum accuracy is considered, as described earlier when presenting the results of Figure~\ref{fig:couetteConvModes}. The results show that, when enough modes are considered within the PGD framework, the error $\varepsilon_{\texttt{PGD}}$ converges with near an optimal rate $h^{\nDeg+1}$. 

%______________________________________________________________________
\subsection{Stokes flow around two circular cylinders}

\begin{figure}[!tb]
	\centering
	\includegraphics[width=0.9\textwidth]{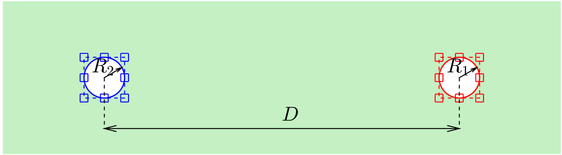}
	\caption{Detail of the computational domain for the solution of the Stokes flow around two circular cylinders.}
	\label{fig:swimmerSetup}
\end{figure}

The second example is inspired by the studies of~\cite{stone1996propulsion,alouges2009optimal} on the analysis of micro-swimmers. The problem involves the computation of the Stokes flow in a rectangular channel of dimension $40 \times 14$ with two circular cylinders of radius $R_1$ and $R_2$, where the distance between their centres is $D$, as represented in Figure~\ref{fig:swimmerSetup}. Slip boundary conditions are considered on the surface of the cylinders, a free slip boundary condition on the top and bottom boundaries, a imposed horizontal velocity of unit magnitude on the left boundary and a homogeneous Neumann boundary condition on the right part of the boundary.

The geometric parameters considered are the radius of the cylinders and the distance between their centres. In addition, as done in the context of microswimmers, it is assumed that the variation of the two radius is linked so that $R_1^2 + R_2^2$ remains constant, that is the area enclosed by the two cylinders is constant. It is worth noting that the variation of these geometric parameters involve the variation of all the control points (i.e. 32 parameters) of the two NURBS describing the circles represented in Figure~\ref{fig:swimmerSetup}. However, as done in the previous example, it is possible to re-parametrise the motion in terms of only two parameters, $\mu_1 \in  \I_1 = [-1,1]$ which controls the radius of the cylinders and $\mu_2 \in  \I_2 = [-1,1]$ controlling the distance between the cylinders.

The reference configuration corresponds to $R_1=R_2=0.8$ ($\mu_1=0$) and $D=14$ ($\mu_2=0$). The minimum and maximum values for the radius of the cylinders are 0.3578 ($\mu_1=-1$) and 1.0733 ($\mu_1=1$) which, given the link between $R_1$ and $R_2$, correspond to the cases where the area of one of the circles is 90\% and 10\% respectively of the total area occupied by both circles. For the distance, the maximum and minimum values are 15.5 ($\mu_2=-1$) and 12.5 ($\mu_2=1$) respectively.

\begin{figure}[!tb]
	\centering
	\subfigure[$\mu_1=\mu_2=0$]{\includegraphics[width=\textwidth]{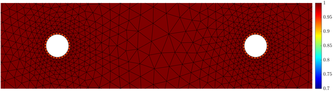}}  \\
	\subfigure[$\mu_1=\mu_2=-1$]{\includegraphics[width=\textwidth]{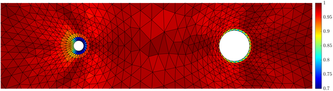}} \\
	\subfigure[$\mu_1=\mu_2=1$]{\includegraphics[width=\textwidth]{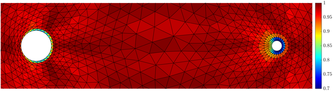}}
	\caption{Computational mesh and quality for three values of the geometric parameters $\mu_1$ and $\mu_2$.}
	\label{fig:swimmerMeshes}
\end{figure}

A computational mesh with 2,338 triangular elements is generated for the reference configuration. For a cubic degree of approximation (for geometry and velocity), Figure~\ref{fig:swimmerMeshes} shows the quality of the mesh for the reference configuration and for two deformed configurations corresponding to the extreme cases with $\mu_1=\mu_2=-1$ and $\mu_1=\mu_2=1$.
It can be observed that, even for the large deformations corresponding to the extreme cases of Figures~\ref{fig:swimmerMeshes} (b) and (c), the minimum quality (measured as the scaled Jacobian) is near 0.7, which is similar to the quality observed in the previous example for a simpler problem.

Next, the separation of the matrix $\mat{H}_{\bmu}$, defined in Equation~\eqref{eq:H}, is considered. As discussed previously $\mat{H}_{\bmu}$ does not generally admit an exact separable expression and therefore a separable approximation is computed here via the higher-order PGD-projection described in \cite{DM-MZH:15}. 

\begin{figure}[!tb]
	\centering
	\subfigure[$m=1$]{\includegraphics[width=0.32\textwidth]{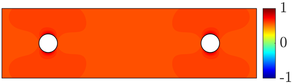}}
	\subfigure[$m=2$]{\includegraphics[width=0.32\textwidth]{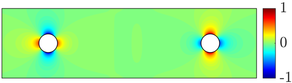}}
	\subfigure[$m=3$]{\includegraphics[width=0.32\textwidth]{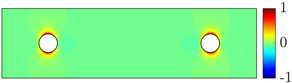}}
	\subfigure[$m=4$]{\includegraphics[width=0.32\textwidth]{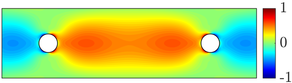}}
	\subfigure[$m=5$]{\includegraphics[width=0.32\textwidth]{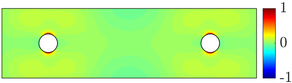}}
	\subfigure[$m=6$]{\includegraphics[width=0.32\textwidth]{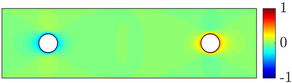}}
	\caption{First six normalised spatial modes of $[\mat{H}_{\bmu}]_{11}$ on the mesh of Figure~\ref{fig:swimmerMeshes} (a) with 25 elements in $\I_1$, 45 elements in $\I_2$ and with $\nDeg=4$.}
	\label{fig:swimmer_H11_Modes}
\end{figure}

\begin{figure}[!tb]
	\centering
	\subfigure[$m=1$]{\includegraphics[width=0.32\textwidth]{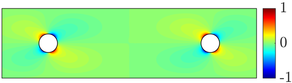}}
	\subfigure[$m=2$]{\includegraphics[width=0.32\textwidth]{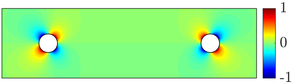}}
	\subfigure[$m=3$]{\includegraphics[width=0.32\textwidth]{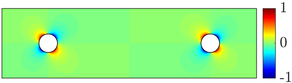}}
	\subfigure[$m=4$]{\includegraphics[width=0.32\textwidth]{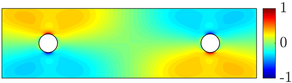}}
	\subfigure[$m=5$]{\includegraphics[width=0.32\textwidth]{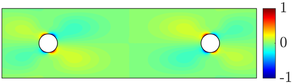}}
	\subfigure[$m=6$]{\includegraphics[width=0.32\textwidth]{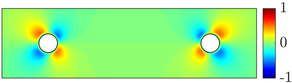}}
	\caption{First six normalised spatial modes of $[\mat{H}_{\bmu}]_{12}$ on the mesh of Figure~\ref{fig:swimmerMeshes} (a) with 25 elements in $\I_1$, 45 elements in $\I_2$ and with $\nDeg=4$.}
	\label{fig:swimmer_H12_Modes}
\end{figure}

\begin{figure}[!tb]
	\centering
	\subfigure[$m=1$]{\includegraphics[width=0.32\textwidth]{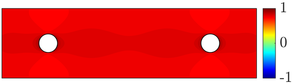}}
	\subfigure[$m=2$]{\includegraphics[width=0.32\textwidth]{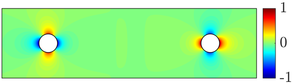}}
	\subfigure[$m=3$]{\includegraphics[width=0.32\textwidth]{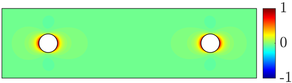}}
	\subfigure[$m=4$]{\includegraphics[width=0.32\textwidth]{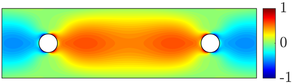}}
	\subfigure[$m=5$]{\includegraphics[width=0.32\textwidth]{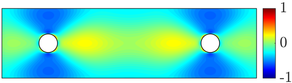}}
	\subfigure[$m=6$]{\includegraphics[width=0.32\textwidth]{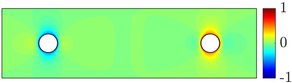}}
	\caption{First six normalised spatial modes of $[\mat{H}_{\bmu}]_{22}$ on the mesh of Figure~\ref{fig:swimmerMeshes} (a) with 25 elements in $\I_1$, 45 elements in $\I_2$ and with $\nDeg=4$.}
	\label{fig:swimmer_H22_Modes}
\end{figure}

Figure~\ref{fig:swimmer_H11_Modes} shows the first six normalised spatial modes of the component $[\mat{H}_{\bmu}]_{11}$.
Similarly, Figures~\ref{fig:swimmer_H12_Modes} and \ref{fig:swimmer_H22_Modes} show the first six normalised spatial modes of the component $[\mat{H}_{\bmu}]_{12}$ and $[\mat{H}_{\bmu}]_{22}$ respectively.

It is worth noting that some of the modes in Figure~\ref{fig:swimmer_H11_Modes} resemble the modes obtained in the previous example and represented in Figure~\ref{fig:couetteMesh_H3P4_H11_Mode}, whereas other modes display a completely different spatial variation. This indicates that the similar modes are the ones that carry information about the change of the radius of the cylinders whereas the different ones are related to the variation of the distance between the cylinders. 

Contrary to the previous example, it is apparent that more modes are necessary to describe the global behaviour of the matrix $\mat{H}_{\bmu}$. For instance, in the previous example only two modes were needed to describe the global information of the matrix $\mat{H}_{\bmu}$ whereas now the fourth mode of all the components of $\mat{H}_{\bmu}$ and the fifth mode of $[\mat{H}_{\bmu}]_{22}$ contain global information. Also, contrary to the previous example, it is interesting to observe the different behaviour of the modes associated to the diagonal terms of the matrix $\mat{H}_{\bmu}$. This is due to the more complex motion induced by the geometric parameters, compared to the previous example. 

\begin{figure}[!tb]
	\centering
	\subfigure[$\mu_1$]{\includegraphics[width=0.49\textwidth]{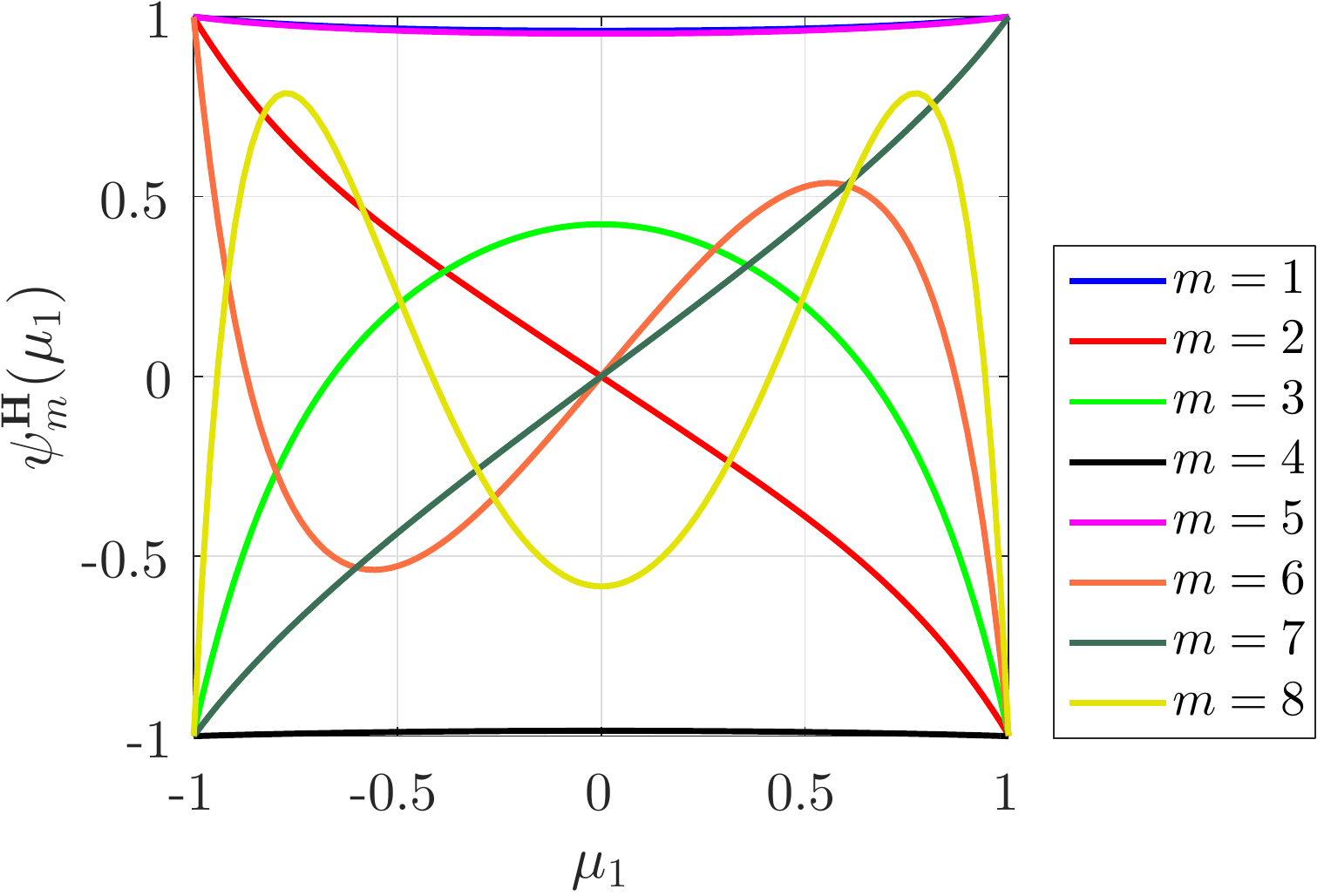}}
	\subfigure[$\mu_2$]{\includegraphics[width=0.49\textwidth]{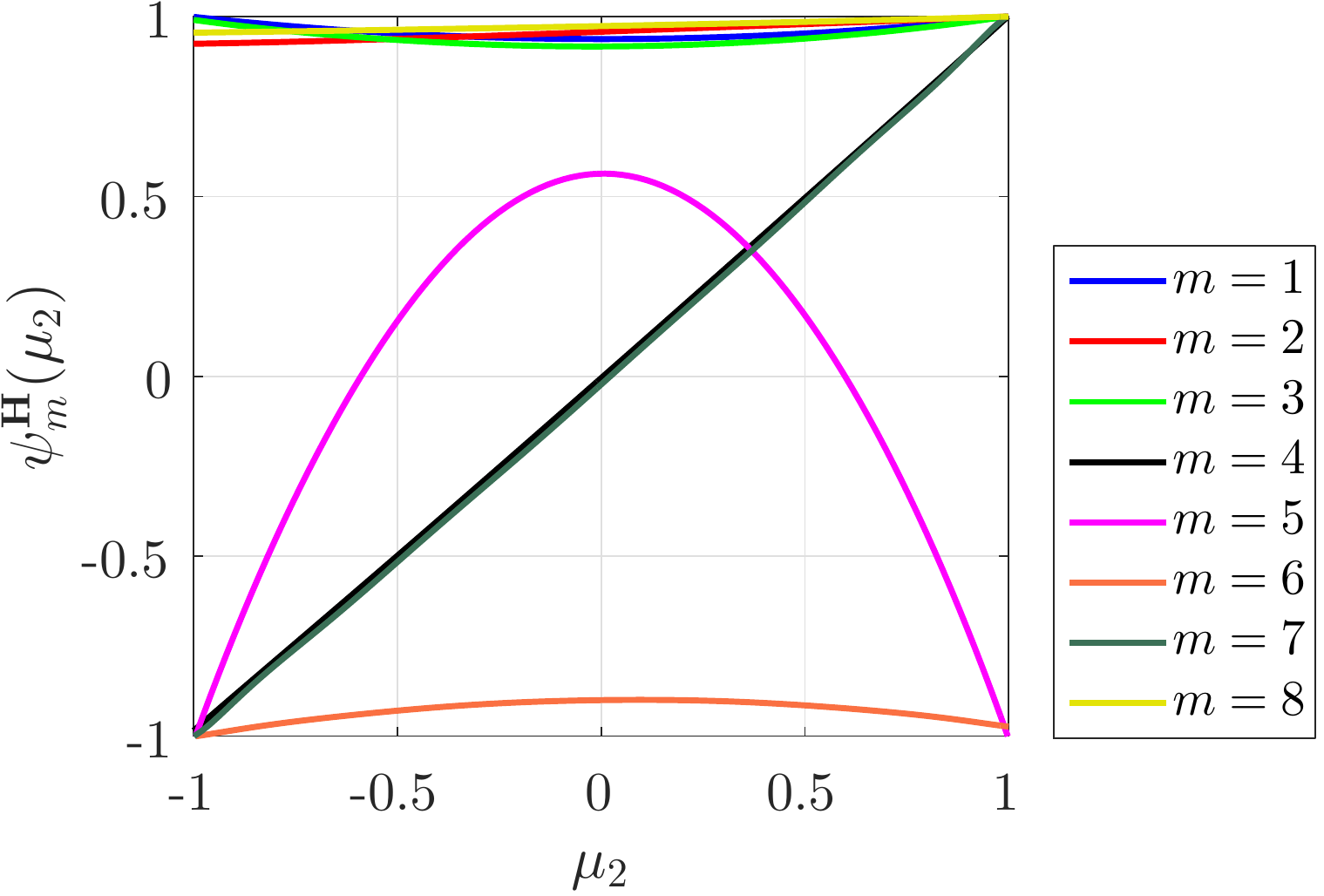}}
	\caption{First eight parametric modes of $\mat{H}_{\bmu}$ on the mesh of Figure~\ref{fig:swimmerMeshes} (a) with 25 elements in $\I_1$, 45 elements in $\I_2$ and with $\nDeg=4$.}
	\label{fig:swimmer_HXX_paramModes}
\end{figure}

The first eight normalised parametric modes of $\mat{H}_{\bmu}$ are represented in Figure~\ref{fig:swimmer_HXX_paramModes}. 
It is worth recalling that the same parametric modes are associated to all the components of the matrix $\mat{H}_{\bmu}$. It is interesting to observe the similar qualitative behaviour of the parametric modes of $\mu_1$, represented in Figure~\ref{fig:swimmer_HXX_paramModes} (a), and the parametric modes obtained in the previous example, represented in Figure~\ref{fig:couetteMesh_H3P4_HXX_ParamModes}. This is expected as in both cases these modes are related to the variation of the radius of a circle. The parametric modes associated to the distance between the cylinders, represented in Figure~\ref{fig:swimmer_HXX_paramModes} (b), display a less oscillatory character than the modes associated to the variation of the radius of a circle.

\begin{figure}[!tb]
	\centering
	\includegraphics[width=0.45\textwidth]{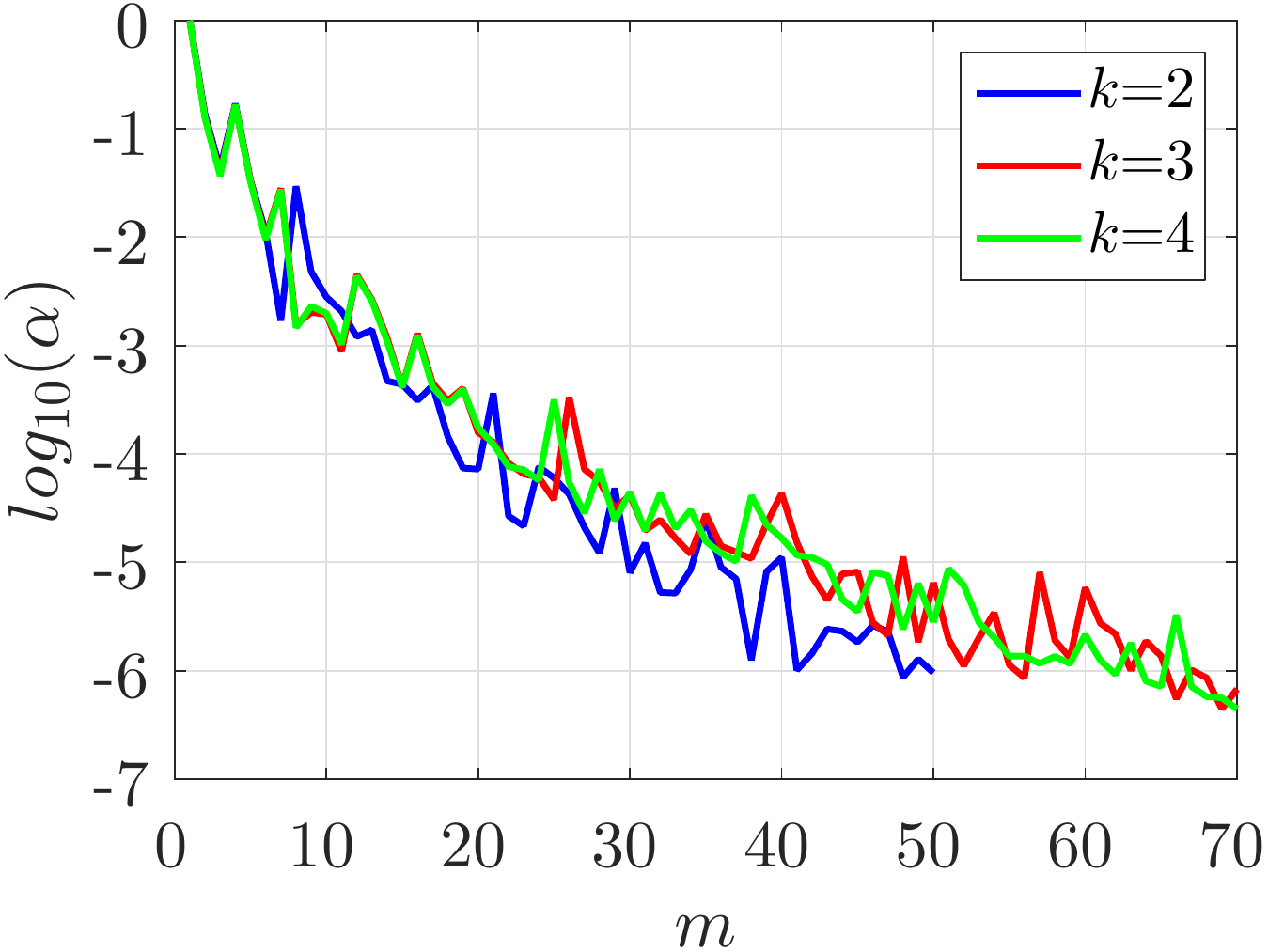}
	\caption{Evolution of the amplitude of the spatial modes $\alpha_m$ of the matrix $\mat{H}_{\bmu}$ using different degrees of approximation.}
	\label{fig:swimmer_AmplitudeH11}
\end{figure}

To illustrate the increased complexity due to the introduction of two geometric parameters, Figure~\ref{fig:swimmer_AmplitudeH11} shows the amplitude, $\alpha_m$, corresponding to the mode $m$ of the separation of $\mat{H}_{\bmu}$, computed as the product of the Euclidean norms of the spatial and parametric functions.

The amplitude, $\alpha_m$, corresponding to the mode $m$ of the separation of $\mat{H}_{\bmu}$, is computed as the product of the Euclidean norms of the spatial and parametric functions. Figure~\ref{fig:couette_AmplitudeH11} shows the amplitudes of $[\mat{H}_{\bmu}]_{11}$ using four different meshes and three different degrees of approximation.
The results show that, for this example, the number of terms required to obtain a separable approximation of the matrix $\mat{H}_{\bmu}$ using the higher-order PGD-projection is significantly higher than in the previous example. To provide a decrease in the amplitude $\alpha_m$ of six orders of magnitude the number of modes vary from 50 and 63 for an order of approximation ranging from $k=2$ to $k=4$, showing again that the number of modes is not highly dependent upon the spatial discretisation. 

\begin{figure}[!tbp]
	\centering
	\subfigure[Velocity, $m=1$]{\includegraphics[width=0.44\textwidth]{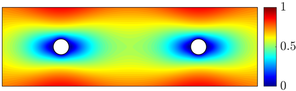}}
	\subfigure[Pressure, $m=1$]{\includegraphics[width=0.44\textwidth]{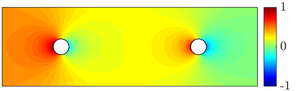}}
	\subfigure[Velocity, $m=2$]{\includegraphics[width=0.44\textwidth]{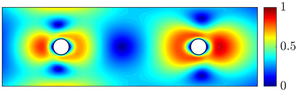}}
	\subfigure[Pressure, $m=2$]{\includegraphics[width=0.44\textwidth]{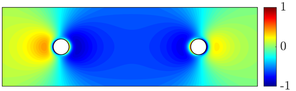}}
	\subfigure[Velocity, $m=3$]{\includegraphics[width=0.44\textwidth]{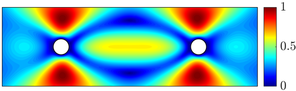}}
	\subfigure[Pressure, $m=3$]{\includegraphics[width=0.44\textwidth]{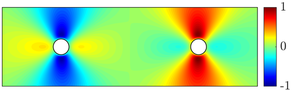}}
	\subfigure[Velocity, $m=4$]{\includegraphics[width=0.44\textwidth]{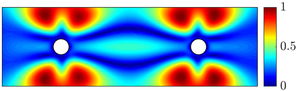}}
	\subfigure[Pressure, $m=4$]{\includegraphics[width=0.44\textwidth]{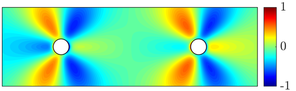}}
	\subfigure[Velocity, $m=5$]{\includegraphics[width=0.44\textwidth]{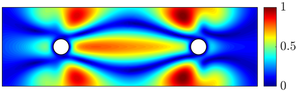}}
	\subfigure[Pressure, $m=5$]{\includegraphics[width=0.44\textwidth]{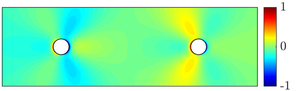}}
	\subfigure[Velocity, $m=6$]{\includegraphics[width=0.44\textwidth]{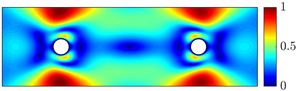}}
	\subfigure[Pressure, $m=6$]{\includegraphics[width=0.44\textwidth]{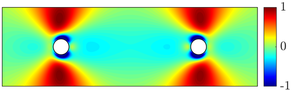}}
	\caption{First six normalised spatial modes of the magnitude of the velocity (left) and pressure (right) computed on the mesh of Figure~\ref{fig:swimmerMeshes} (a) with 25 elements in $\I_1$, 45 elements in $\I_2$ and with $\nDeg=4$.}
	\label{fig:swimmerModes}
\end{figure}

Using the separation of the matrix $\mat{H}_{\bmu}$, the Stokes flow around the two cylinders is computed with the proposed approach to obtain the generalised solution. The first six normalised spatial modes of the magnitude of the velocity and the pressure fields are shown in Figure~\ref{fig:swimmerModes}.
The simulation was performed using the mesh of Figure~\ref{fig:swimmerMeshes} (a), 25 equally-spaced elements in $\I_1 = [-1,1]$, 45 equally-spaced elements in $\I_2 = [-1,1]$ and a degree of approximation $\nDeg=4$. 

\begin{figure}[!tb]
	\centering
	\subfigure[$\mu_1$]{\includegraphics[width=0.49\textwidth]{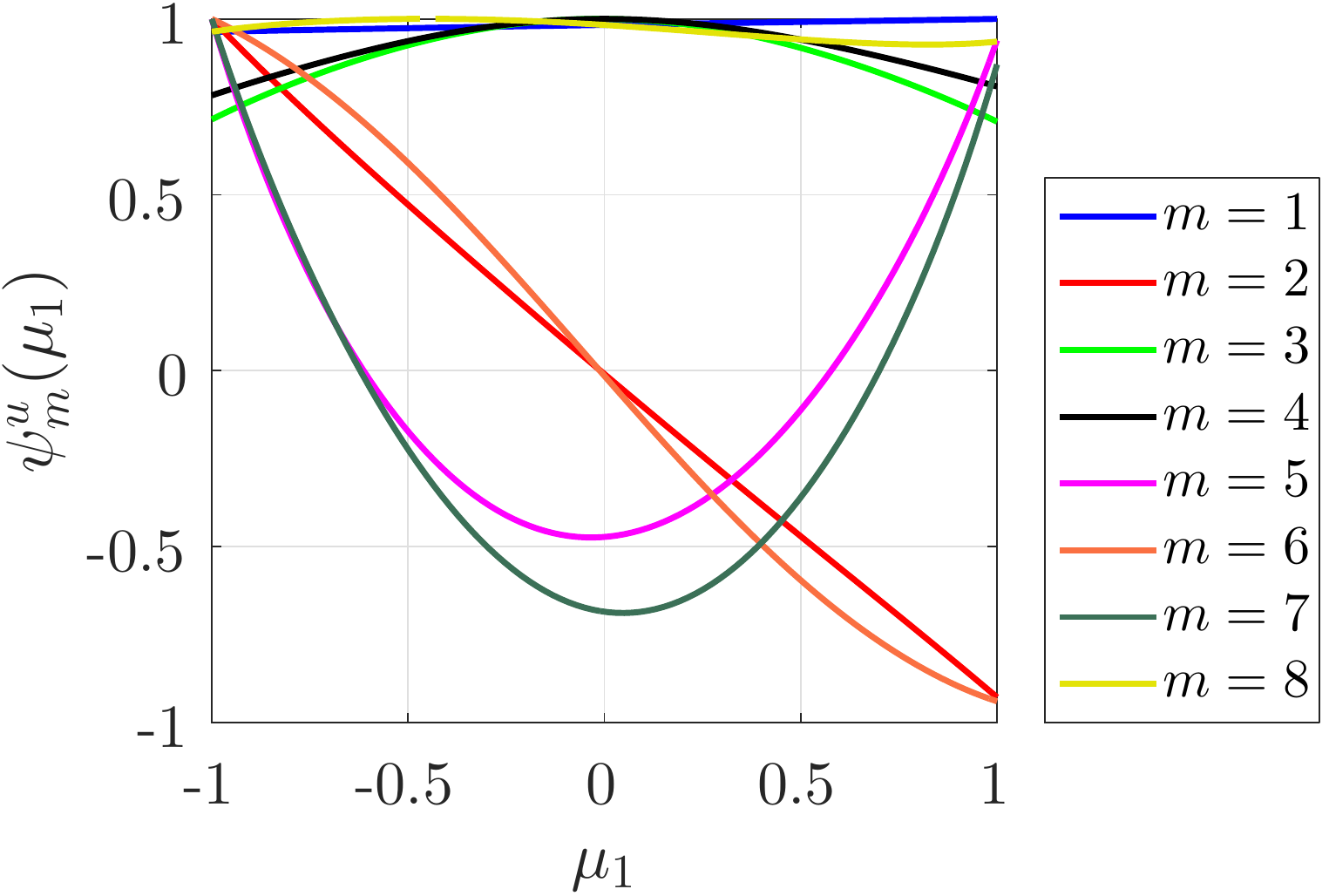}}
	\subfigure[$\mu_2$]{\includegraphics[width=0.49\textwidth]{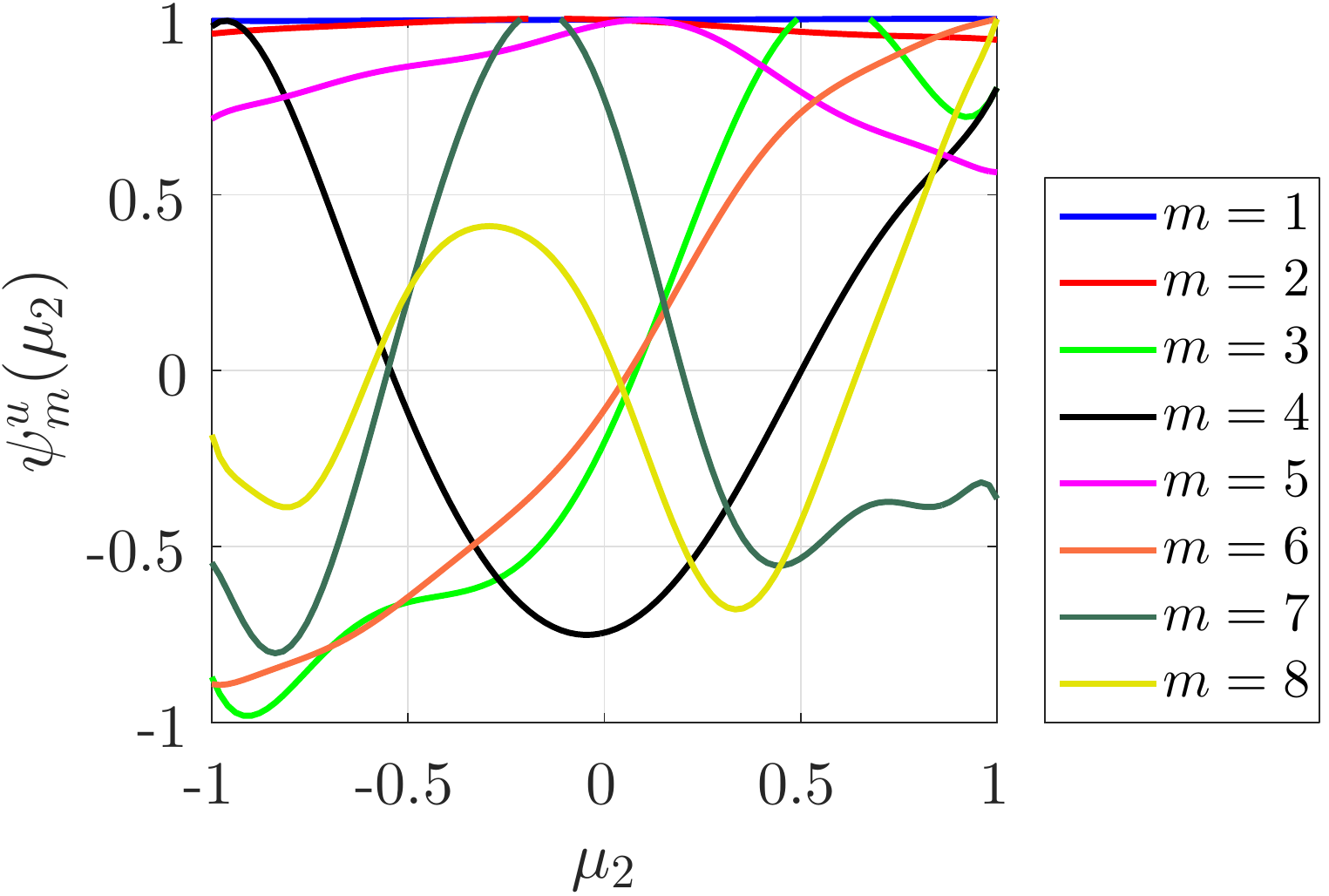}}
	\caption{First eight parametric modes of the Stokes flow around two circular cylinders on the mesh of Figure~\ref{fig:swimmerMeshes} (a) with 25 elements in $\I_1$, 45 elements in $\I_2$ and with $\nDeg=4$.}
	\label{fig:swimmer_paramModes}
\end{figure}

The first eight normalised parametric modes associated to the spatial modes of Figure~\ref{fig:swimmerModes} are represented in Figure~\ref{fig:swimmer_paramModes}. 
It is worth recalling that the same parametric mode is associated to all the components of the velocity and the pressure fields. 

Contrary to the previous example, it can be observed that the spatial modes for the velocity and the pressure do not resemble the spatial modes of the matrix $\mat{H}_{\bmu}$, illustrating the increased complexity of the current problem. In addition, it is worth noting that the parametric modes associated to the first parameter $\mu_1$ shows a more oscillatory character than $\mu_2$ when the separation of the matrix $\mat{H}_{\bmu}$ is studied (as shown in Figure~\ref{fig:swimmer_HXX_paramModes}), whereas the the second parameter $\mu_2$ shows a more oscillatory character than $\mu_1$ when the separation of the velocity and pressure is considered (as shown in Figure~\ref{fig:swimmer_paramModes}). This indicates that the flow around the two cylinders changes slowly when the radius of the cylinders is varied whereas the flow changes more rapidly when the distance between the cylinders is varied.

\begin{figure}[!tb]
	\centering
	\subfigure[Velocity, $\mu_1=\mu_2=0$]{\includegraphics[width=0.49\textwidth]{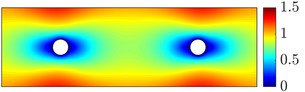}}
	\subfigure[Pressure, $\mu_1=\mu_2=0$]{\includegraphics[width=0.49\textwidth]{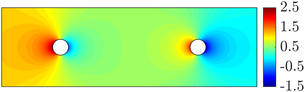}}
	\subfigure[Velocity, $\mu_1=\mu_2=-1$]{\includegraphics[width=0.49\textwidth]{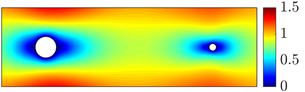}}
	\subfigure[Pressure, $\mu_1=\mu_2=-1$]{\includegraphics[width=0.49\textwidth]{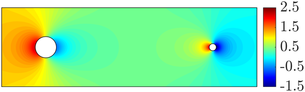}}
    \subfigure[Velocity, $\mu_1=\mu_2=1$]{\includegraphics[width=0.49\textwidth]{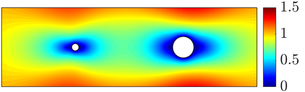}}
	\subfigure[Pressure, $\mu_1=\mu_2=1$]{\includegraphics[width=0.49\textwidth]{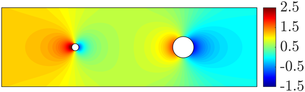}}
	\caption{Magnitude of the velocity (left) and pressure (right) computed on the three configurations shown in Figure~\ref{fig:swimmerMeshes} with 25 elements in $\I_1$, 45 elements in $\I_2$ and with $\nDeg=4$.}
	\label{fig:swimmerSolutions}
\end{figure}

Next, the solutions obtained with the proposed PGD framework are represented for different values of the geometric parameters. Figure~\ref{fig:swimmerSolutions} shows the magnitude of the velocity and the pressure fields for the three configurations shown in Figure~\ref{fig:swimmerMeshes}.

\begin{figure}[!tb]
	\centering
	\subfigure[Velocity, $\mu_1=\mu_2=0$]{\includegraphics[width=0.49\textwidth]{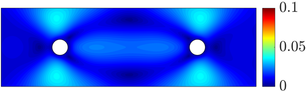}}
	\subfigure[Pressure, $\mu_1=\mu_2=0$]{\includegraphics[width=0.49\textwidth]{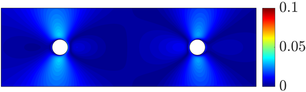}}
	\subfigure[Velocity, $\mu_1=\mu_2=-1$]{\includegraphics[width=0.49\textwidth]{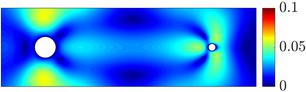}}
	\subfigure[Pressure, $\mu_1=\mu_2=-1$]{\includegraphics[width=0.49\textwidth]{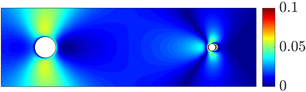}}
	\subfigure[Velocity, $\mu_1=\mu_2=1$]{\includegraphics[width=0.49\textwidth]{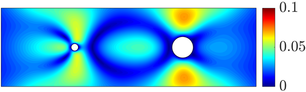}}
	\subfigure[Pressure, $\mu_1=\mu_2=1$]{\includegraphics[width=0.49\textwidth]{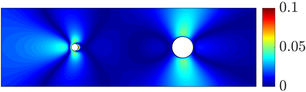}}
	\caption{Difference between the standard FE solution and the PGD solution for the magnitude of the velocity (left) and pressure (right) computed on the three configurations shown in Figure~\ref{fig:swimmerSolutions}.}
	\label{fig:swimmerErrMaps}
\end{figure}

To illustrate the accuracy of the proposed approach, Figure~\ref{fig:swimmerErrMaps} shows the difference between the standard FE solution and the PGD solution for the magnitude of the velocity and the pressure corresponding to the three configurations shown in Figure~\ref{fig:swimmerMeshes}.
For all the configurations, the results obtained with the proposed PGD framework are in excellent agreement with the results obtained using a standard FE solver. For the first configuration displayed in Figure~\ref{fig:swimmerErrMaps}, corresponding to $\mu_1=\mu_2=0$, the difference between the magnitude of the velocity field using FEM and the PGD in the $\eltwo(\Omega^{\bmu})$ norm is 0.0143 and the difference between the pressure field using FEM and the PGD in the $\eltwo(\Omega^{\bmu})$ norm is 0.0277. For the second configuration of Figure~\ref{fig:swimmerErrMaps}, corresponding to $\mu_1=\mu_2=-1$, the difference between the magnitude of the velocity and pressure fields using FEM and the PGD in the $\eltwo(\Omega^{\bmu})$ norm are 0.0267 and 0.0592 respectively. Finally, for the last configuration of Figure~\ref{fig:swimmerErrMaps}, corresponding to $\mu_1=\mu_2=1$, the difference between the magnitude of the velocity and pressure fields using FEM and the PGD in the $\eltwo(\Omega^{\bmu})$ norm are 0.0218 and 0.0787 respectively. It is worth emphasising that the accurate results obtained, with differences between PGD and standard FE simulations lower than 8\% for the whole range of geometric configurations, have been obtained with very coarse meshes in the parametric spaces, only 25 elements in $\I_1$ and 45 elements in $\I_2$.

\subsection{Stokes flow around an arbitrarily shaped vesicle}

The last example is inspired by the study of vesicles suspended in a viscous flow presented in~\cite{veerapaneni2011fast}. The characterisation of such flows is of interest in many biomechanical applications and the simulations often require the computation to be performed for a large variety of geometric configurations. This example is used to demonstrate the applicability of the proposed technique in three dimensions by using three geometric parameters that lead to substantial variations in the geometric model.

The problem considered here involves the simulation of the Stokes flow around an arbitrarily shaped vesicle in a channel of dimension $10 \times 10 \times 20$.  Slip boundary conditions are considered on the surface of the vesicle, an imposed vertical velocity of unit magnitude in the inflow part of the boundary ($z=-10$), a homogeneous Neumann boundary condition on the outflow part of the boundary ($z=10$) and free slip boundary conditions on the remaining parts of the boundary. 

\begin{figure}[!tb]
	\centering
	\includegraphics[width=0.55\textwidth]{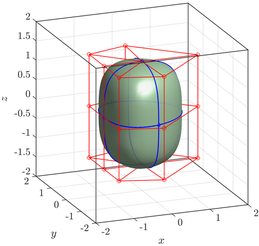}
	\caption{Geometric model of a vesicle, showing the control net and the breaklines.}
	\label{fig:vesicleNURBSref}
\end{figure}

The generic vesicle considered is modelled using a single degenerate cubic NURBS with 20 control points and four patches, as illustrated in Figure~\ref{fig:vesicleNURBSref}. 
The control points are arranged in three groups, according to their $z$ coordinate being positive, zero or negative. The geometric parametrisation considers the motion of all the three groups independently by using three parameters, $\mu_1$, $\mu_2$ and $\mu_3$ respectively. Each parameter is used to vary the radial position of the control points whilst maintaining its $z$ coordinate. For the control points with positive $z$ coordinate the radial coordinate is given by $r + 3\mu_1/4$ with $\mu_1 \in [-1,1]$. Similarly for the control points with zero or positive $z$ coordinate the radial coordinate is $r + 3\mu_2/4$ and $r + 3\mu_3/4$ respectively, with $\mu_2 \in [-1,1]$ and $\mu_3 \in [-1,1]$.

The reference configuration, corresponding to $\bmu=(\mu_1,\mu_2,\mu_3) = (0,0,0)$, is shown in Figure~\ref{fig:vesicleNURBSref}.

The generalised PGD solution is computed using a tetrahedral mesh with 6,712 elements, 10 equally-spaced elements in each parametric interval $\I_1 = \I_2 = \I_3 = [-1,1]$ and a degree of approximation $\nDeg=4$. 

\begin{figure}[!tb]
	\centering
	\subfigure[$\bmu=(-1, -1,  1)$]{\includegraphics[width=0.32\textwidth]{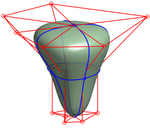}}
	\subfigure[$\bmu=(-1,  1,  1)$]{\includegraphics[width=0.32\textwidth]{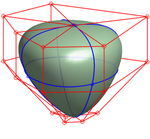}}
	\subfigure[$\bmu=( 0, -1,  1)$]{\includegraphics[width=0.32\textwidth]{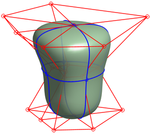}}
	\subfigure[$\bmu=( 0,  1,  0)$]{\includegraphics[width=0.32\textwidth]{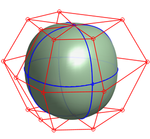}}
	\subfigure[$\bmu=( 1,  0, -1)$]{\includegraphics[width=0.32\textwidth]{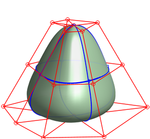}}
	\subfigure[$\bmu=( 1,  0,  0)$]{\includegraphics[width=0.32\textwidth]{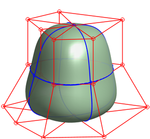}}
	\caption{Geometric model of six different vesicle obtained by deforming the reference configuration of Figure~\ref{fig:vesicleNURBSref}.}
	\label{fig:vesicleNURBS}
\end{figure}

\begin{figure}[!tb]
	\centering
	\subfigure[$\bmu=(-1, -1,  1)$]{\includegraphics[width=0.32\textwidth]{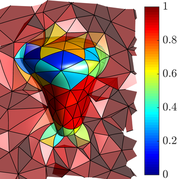}}
	\subfigure[$\bmu=(-1,  1,  1)$]{\includegraphics[width=0.32\textwidth]{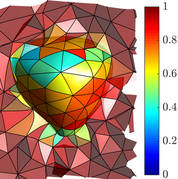}}
	\subfigure[$\bmu=( 0, -1,  1)$]{\includegraphics[width=0.32\textwidth]{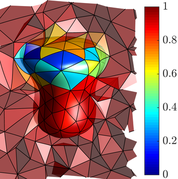}}
	\subfigure[$\bmu=( 0,  1,  0)$]{\includegraphics[width=0.32\textwidth]{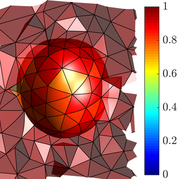}}
	\subfigure[$\bmu=( 1,  0, -1)$]{\includegraphics[width=0.32\textwidth]{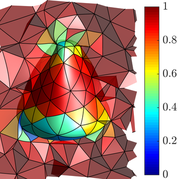}}
	\subfigure[$\bmu=( 1,  0,  0)$]{\includegraphics[width=0.32\textwidth]{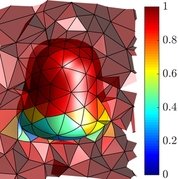}}
	\caption{Detail of the six deformed high-order meshes corresponding to the configurations shown in Figure~\ref{fig:vesicleNURBS}.}
	\label{fig:vesicleMesh}
\end{figure}

To illustrate the variation in the geometry induced by the selected geometric parameters, Figure~\ref{fig:vesicleNURBS} shows six different geometric configurations using different values for the geometric parameters.
It can be observed that the variations of the control points considered lead to substantial changes in the geometric model. These variations induce a deformation of the mesh generated for the reference configuration, computed in a separated form by employing the strategy described in Section~\ref{sbc:domainDisplacement}. Figure~\ref{fig:vesicleMesh} shows six deformed meshes near the vesicle that correspond to the six variations of the geometric parameters depicted in Figure~\ref{fig:vesicleNURBS}.
The colours represent the quality of the elements, measured as the scaled Jacobian.

Following the strategy described in this work, the matrix $\mat{H}_{\bmu}$ is separated first using the higher-order PGD-projection. To provide a decrease in the amplitude $\alpha_m$ of five orders of magnitude, 50 modes are required in this example. It is worth noting that this is almost the same number of modes required for the previous two dimensional example to provide the same degree in the amplitude of the modes of the matrix $\mat{H}_{\bmu}$.

Using the separation of the matrix $\mat{H}_{\bmu}$, the Stokes flow around the vesicle is computed with the proposed approach to obtain the generalised solution. A total of 160 modes are computed, being the amplitude of the last mode three orders of magnitude lower than the amplitude of the first mode. 

\begin{figure}[!tb]
	\centering
	\subfigure[$\bmu=(-1, -1,  1)$]{\includegraphics[width=0.32\textwidth]{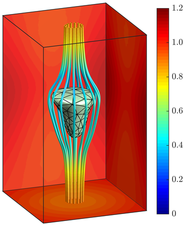}}
	\subfigure[$\bmu=(-1,  1,  1)$]{\includegraphics[width=0.32\textwidth]{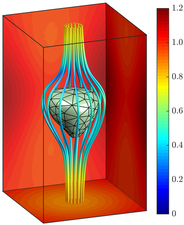}}
	\subfigure[$\bmu=( 0, -1,  1)$]{\includegraphics[width=0.32\textwidth]{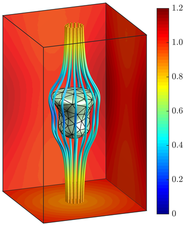}}
	\subfigure[$\bmu=( 0,  1,  0)$]{\includegraphics[width=0.32\textwidth]{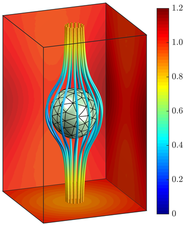}}
	\subfigure[$\bmu=( 1,  0, -1)$]{\includegraphics[width=0.32\textwidth]{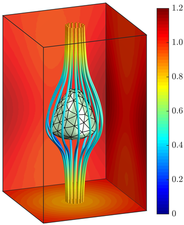}}
	\subfigure[$\bmu=( 1,  0,  0)$]{\includegraphics[width=0.32\textwidth]{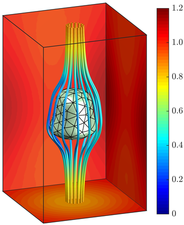}}
	\caption{Magnitude of the velocity and isolines for the six configurations shown in Figure~\ref{fig:vesicleNURBS}.}
	\label{fig:vesicleVelo}
\end{figure}

\begin{figure}[!tb]
	\centering
	\subfigure[$\bmu=(-1, -1,  1)$]{\includegraphics[width=0.32\textwidth]{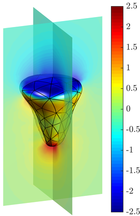}}
	\subfigure[$\bmu=(-1,  1,  1)$]{\includegraphics[width=0.32\textwidth]{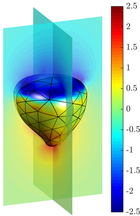}}
	\subfigure[$\bmu=( 0, -1,  1)$]{\includegraphics[width=0.32\textwidth]{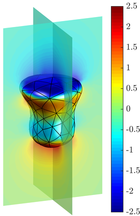}}
	\subfigure[$\bmu=( 0,  1,  0)$]{\includegraphics[width=0.32\textwidth]{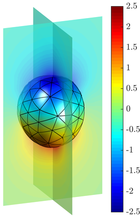}}
	\subfigure[$\bmu=( 1,  0, -1)$]{\includegraphics[width=0.32\textwidth]{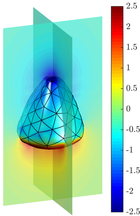}}
	\subfigure[$\bmu=( 1,  0,  0)$]{\includegraphics[width=0.32\textwidth]{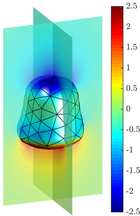}}
	\caption{Pressure field for the six configurations shown in Figure~\ref{fig:vesicleNURBS}.}
	\label{fig:vesiclePress}
\end{figure}

Figures~\ref{fig:vesicleVelo} and \ref{fig:vesiclePress} show the velocity and pressure fields for the six geometric configurations shown in Figure~\ref{fig:vesicleNURBS}. 
These configurations are obtained in the on-line phase as particularisations of the genearlised, six dimensional, PGD solution. Both the velocity and pressure fields illustrate the ability of the proposed PGD approach to capture significant changes in the flow field induced by geometric variations of the CAD model. 

To quantify the accuracy of the proposed technique in a more complex problem in three dimensions and with three geometric parameters, the particularised solutions for the velocity and pressure are compared to standard finite element computations performed on the deformed configurations for the six cases displayed in Figures~\ref{fig:vesicleVelo} and \ref{fig:vesiclePress}. For the velocity field, the difference between the PGD and the finite element solution, measured in the $\eltwo(\Omega)$ norm, is 0.7\%, 0.8\%, 0.8\%, 0.4\%, 0.6\% and 0.7\% for the six cases shown in Figure~\ref{fig:vesicleVelo} respectively. Similarly, for the pressure field, the difference between the PGD and the finite element solution, measured in the $\eltwo(\Omega)$ norm, is 7\%, 7\%, 7\%, 3\%, 5\% and 6\% for the six cases shown in Figure~\ref{fig:vesiclePress} respectively. It is worth emphasising that the higher accuracy observed in the velocity is due to the use of a higher interpolation degree for the velocity field, compared to the pressure field, in order to satisfy the LBB conditions.

%______________________________________________________________________
\section{Concluding remarks}

A computational framework for the computation of off-line solutions for a set of parameters describing the geometry of a domain has been presented. The proposed approach considers as parameters of the generalised PGD solution the control points of the CAD boundary representation of the computational domain. A mapping between a reference configuration and the current configuration is proposed by interpreting the geometric changes as a displacement field derived from a linear elastic problem. The key aspect of the proposed approach is that the displacement field is explicitly written in a separable form. This approach enables the incorporation of the PGD rationale in a high-order finite element context.

The potential of the proposed approach is shown for a variety of problems involving the solution of the Stokes equation in geometrically parametrised domains, both in two and three dimensions. The problems, of increasing difficulty show the optimal approximation properties of the method and its ability to accurately capture the flow features for significant changes of the geometric model. For the most complex problem considered in this work, the generalised solution computed with three geometric parameters in three dimensions show good agreement when compared to a standard finite element computation, with errors in the velocity field lower than 1\% and errors in the pressure field below 8\%.

\appendix
%______________________________________________________________________
\section{Integration within a CAD environment in 3D}  \label{app:cadIntegration3D}

In three dimensions, the boundary of the parameterized domain, $\partial \Omega^{\bmu}$, is assumed to be described by using a set NURBS surfaces $\{\surface[j]^\bmu\}_{j=1,\dotsc, M}$, being $M$ the total number of surfaces, namely
\begin{equation*}
	\partial \Omega^{\bmu} =  \bigcup_{j=1}^{M} \surface[j]^\bmu([0,1]^2).
\end{equation*}

Next, the necessary concepts about NURBS surfaces are briefly recalled and the minimum changes that are required to extend the technique presented in Section \ref{sc:cadIntegration} to three dimensional domains are detailed.

%=======================================================================
\subsection{NURBS surfaces} \label{sbc:NURBSsurfaces}
%=======================================================================
	
A NURBS surface of degree $q$ in $\lambda$ and degree $r$ in $\kappa$, is a piecewise rational function defined in parametric form as
\begin{equation*} \label{eq:nurbsSurfaceDefinition}
\surface(\lambda,\kappa) = \sum_{i=0}^{\ncp}\sum_{l=0}^{\mcp}
\bm{B}_{il} R_{i,l}(\lambda,\kappa)
\qquad 0 \leq  \lambda,\kappa \leq 1,
\end{equation*}
where $\{\bm{B}_{il}\}$ are the coordinates of the $(\ncp+1)(\mcp+1)$ \emph{control points} (defining the \emph{control net}) and $\{R_{il}\}$ are rational basis functions defined as
\begin{equation*} \label{eq:rationalBasisSurface}
R_{il}(\lambda,\kappa) =  \nu_{il} \, S_{i,l}^{q,r}(\lambda,\kappa)
\! \biggm/\!\!
\left(\sum_{i=0}^{\ncp}\sum_{l=0}^{\mcp} \nu_{il}\,
S_{i,l}^{q,r}(\lambda,\kappa)\right).
\end{equation*}

In the above expression $\{\nu_{il}\}$ are the control weights associated to the control points and $\{S_{i,l}^{q,r}(\lambda,\kappa)\}$  are the 2D B-spline basis functions of degree $q$ in $\lambda$ and $r$ in $\kappa$. Each 2D B-Spline basis function is defined as a tensor product of 1D basis functions, that is
\begin{equation*} \label{eq:nurbsSurfaceBSpline2D}
S_{i,l}^{q,r}(\lambda,\kappa) :=  C_i^q(\lambda)C_l^r(\kappa).
\end{equation*}

%=======================================================================
\subsection{Geometric parameters}  \label{sbc:geometricParametersSurfaces}
%=======================================================================

The geometric parameters $\bmu\in\bI\subset\mathbb{R}^{\npar}$ are defined as the variations of the original coordinates of the control points of the NURBS surfaces describing the boundary. More precisely, for each NURBS curve $\surface[j]$, with $j=1,\dotsc, M$, having $(\ncp^j + 1)(\mcp^j + 1)$  control points, the undisturbed boundary is characterised by the coordinates of the control points:  $\bm{B}_{il}^j$, for  $i=0,\dotsc , \ncp^j$ and $l=0,\dotsc , \mcp^j$.  The boundary of the distorted domain, $\partial\Omega^\bmu$, is defined by the position of the displaced control points, namely $\bm{B}_{il}^j+\dc_{il}^j$. The displacement range for each control point is characterised by 
\begin{equation*}
\dc_{il}^{il}\in\bI^j_{il} = [\underline{\dx}_{il}^j, \overline{\dx}_{il}^j] \times [\underline{\dy}_{il}^j, \overline{\dy}_{il}^j] \times [\underline{\dz}_{il}^j, \overline{\dz}_{il}^j].
\end{equation*}
In fact, each displacement of a control point on the $j$-th NURBS, $\dc_{il}^j$, might depend upon the parameters and can be written as 
\begin{equation} \label{eq:DeltaControlSurfaces}
\dc_{il}^j=\mu_1^{{il},j}\bm{e}_1 + \mu_2^{{il},j}\bm{e}_2 + \mu_3^{{il},j}\bm{e}_3,
\end{equation}
where $\bm{e}_i$, for $i=1,2,3$, are the unit coordinate vectors. Then $\bmu\in\bI := \bI^1\times\bI^2 \times\dotsb\times\bI^M$, where $\bI^j = \bI^j_1 \times \bI^j_2 \times \cdots \bI^j_{(\ncp^j + 1)(\mcp^j + 1)}$ is the range of variation of the coordinates of the control points of the curve $\surface[j]$.

%=======================================================================
\subsection{Separated representation of the boundary displacement}  \label{sbc:boundaryDisplacementSurfaces}
%=======================================================================

The variation of a control point $\bm{B}_{il}^j$ of a NURBS surface $\surface[j]$, namely $\dc_{il}^j$, changes the definition of the original curve only in the support of the basis function $R_{il}^j$, given by the subspace of the parametric space $[\lambda_i, \lambda_{i+q^j+1}] \times [\kappa_l, \kappa_{l+r^j+1}]$. The modified NURBS surface is parametrised by
\begin{equation*} 
\surface[j]^{\bmu}(\lambda,\kappa) = \sum_{i=0}^{\ncp^j} \sum_{l=0}^{\mcp^j} ( \bm{B}_{il}^j + \dc_{il}^j ) \,R_{il}^j(\lambda,\kappa)
\qquad 0 \leq  \lambda,\kappa \leq 1 .
\end{equation*}

As in the two dimensional case, the displacement of the boundary mesh node $\bX_k=\surface[j](\lambda_k,\kappa_k)$ that belongs to the NURBS curve $\surface[j]$ can be written in separated form as
\begin{equation*}
\delta\bd^j(\bX_k ,\bmu) = \sum_{i=0}^{\ncp^j} \sum_{l=0}^{\mcp^j} \sum_{s=1}^{\nsd} \mu_s^{i,j}\bm{e}_s\,R_{il}^j(\lambda_k,\kappa_k) ,
\end{equation*}
where the dependence of the displacements of the control points in terms of the parameters described in Equation~\eqref{eq:DeltaControlSurfaces} has been used. Moreover, since the NURBS parameterd $(\lambda_k, \kappa_k)$ is only dependent on the spatial coordinates $\bX_k$, and not on the geometric parameters $\bmu$, the previous equation can be written as,
\begin{equation*} 
\delta\bd(\bX_k ,\bmu) = \sum_{j=1}^{M} \delta\bd^j(\bX_k ,\bmu) 
= \sum_{j=1}^{M} \sum_{i=0}^{\ncp^j} \sum_{l=0}^{\mcp^j} \sum_{s=1}^{\nsd} \mu_s^{il,j}\bm{e}_s\,R_{il}^j\bigl(\surface[j]^{-1}(\bX_k)\bigr) ,
\end{equation*}
which characterises the displacement of the boundary nodes and has the desired separated form, as in the two dimensional case given by Equation~\eqref{eq:boundayDisplacementSeparated}.

%______________________________________________________________________
\subsection{Separated representation of the geometric mapping}

The strategy to obtain a separated representation of the geometric mapping is not dependent on the dimensionality of the problem. Therefore, the strategy described in Section~\ref{sbc:domainDisplacement} is also valid in three dimensions. The implementation details are given in \ref{app:solidMechImplementation}.

%______________________________________________________________________
\section{Implementation details of the separated representation of the geometric mapping} \label{app:solidMechImplementation}

To obtain a separated representation of the displacement function in the whole domain, the solid mechanics problem~\eqref{eq:elasticityStrong} is considered. Its discretisation leads to the system of equations~\eqref{eq:solidMechSystemEq}. 

Assuming that the mesh nodes are ordered so that the boundary nodes are first, the vector $\vect{D}$ is given by
\begin{equation*}
\vect{D} =\left( D_1^{X_1}, \dotsc, D_1^{X_{\nsd}}, D_2^{X_1}, \dotsc, D_2^{X_{\nsd}}, \dotsc, D_{|\mathcal{S}|}^{X_1}, \dotsc, D_{|\mathcal{S}|}^{X_{\nsd}} \right)^T ,
\end{equation*}
where $D_k^{X_l}$ is the imposed displacement of node $\bX_k$ in the $X_l$ direction.

As usual in a FE context, the modified system of linear equations to be solved, after accounting for the Dirichlet boundary conditions, is
\begin{equation*} 
\A_{11} \vect{d} = -\A_{12} \vect{D}.
\end{equation*} 

The solution of this system of linear equations can be written as
\begin{equation} \label{eq:nodalDisplacement}
\vect{d} = -\A_{11}^{-1} \A_{12} \vect{D} = \sum_{k=1}^{|\mathcal{S}|} \mat{B}^k \vect{D}_k
\end{equation}
where $\vect{D}_k = \left( D_k^{X_1}, \dotsc, D_k^{X_{\nsd}} \right)^T$ and $\mat{B}^k$ denotes the block of the matrix $\mat{B} := -\A_{11}^{-1} \A_{12}$ containing the columns from $\nsd(k-1)+1$ to $\nsd k$, with dimension $\nsd |\mathcal{S}| \times \nsd$.

Using the separated representation of the imposed displacement \eqref{eq:boundayDisplacementSeparated} for each boundary node, the following separated representation of the nodal values of the displacement is obtained
\begin{equation*} 
\vect{d} = \sum_{k=1}^{|\mathcal{S}|} \mat{B}^k \sum_{j=1}^{M} \sum_{i=0}^{\ncp^j} \bm{\Phi}_i^j(\bmu) \Theta_i^j(\bX_k)  ,
\end{equation*}
leading to the separated representation of the approximation of the displacement function given by Equation~\eqref{eq:separatedDisplacement}.

It is worth noting that the matrix $\mat{B}$ in Equation \eqref{eq:nodalDisplacement} only depends on the spatial discretisation of the original configuration, $\Omega$, and the selected material parameters ($E$ and $\nu$) and it is independent on the geometric parameters $\bmu$. Therefore, it is possible to pre-compute and store the matrix $\mat{B}$ so that the separated representation of the displacement of Equation \eqref{eq:nodalDisplacement} can be computed with a negligible cost for different imposed boundary displacements (i.e. for different configurations $\Omega^{\bmu}$ induced by different variations of the geometric parameters). 

It is important to recall that the dimension of the matrix $\mat{B}$ is $\nsd (\nnodes - |\mathcal{S}| ) \times \nsd |\mathcal{S}|$, which, in practical applications, is much lower than the size of a standard FE matrix, namely $\nsd \nnodes \times \nsd \nnodes$.	

%______________________________________________________________________
\section*{Acknowledgements}

This work is partially supported by the European Union's Horizon 2020 research and innovation programme under the Marie Sk\l odowska--Curie actions (Grant number: 675919) and the Spanish Ministry of Economy and Competitiveness (Grant number: DPI2017-85139-C2-2-R).  The first author also gratefully acknowledges the financial support provided by EPSRC (Grant number: EP/P033997/1). The second and third authors are also grateful for the financial support provided by the Generalitat de Catalunya (Grant number: 2017-SGR-1471).

%________________________________________________________________________
%________________________________________________________________________
%\bibliographystyle{amsplain}
\bibliographystyle{abbrv}
\bibliography{Ref-PGD-Geo}
\end{document}